\documentclass{amsart}
\usepackage{color, graphicx, enumerate, amssymb}
\usepackage{hyperref}
\usepackage{epsfig,wrapfig}
\usepackage{pxfonts}
\usepackage{graphicx}

\usepackage{eucal}

\usepackage[utf8]{inputenc}
\usepackage{ulem}

\usepackage{amsmath}

\usepackage{geometry}
\usepackage{enumitem}

\newtheorem{theorem}{Theorem}[section]

\newtheorem{proposition}[theorem]{Proposition}
\newtheorem{remark}[theorem]{Remark}

\title[]{The role of resonator neuron in the  dynamics of two coupled integrator and resonator neurons of different types of excitability}

\author[M.R. Razvan]{Mohammad Reza Razvan}
\address{Department of Mathematical Sciences, Sharif University of Technology, P.O. Box: 11365-9415, Tehran, Iran}
\email{razvan@sharif.ir}

\author[S. Yasaman]{Somaye Yasaman$^*$}
\address{Department of Mathematical Sciences, Sharif University of Technology, P.O. Box: 11365-9415, Tehran, Iran}
\email{somaye.yasaman@sharif.edu}

\thanks{$^*$Corresponding author: Somaye Yasaman}

\date{\today}

\begin{document}

\begin{abstract}

In this manuscript, a silent resonator neuron is coupled with a spiking integrator neuron through the gap junction, when the coupled neurons are of different types of excitability and none of the coupled neurons exhibit mixed mode oscillations and bursting oscillations intrinsically. By using dynamical systems theory (e.g. the bifurcation theory), all the observed oscillation patterns and the transition mechanisms between them are investigated, when one of the coupling strengths is fixed and the other is varied. It is noticeable that, there is an interval in the parameter space, for the parameter values within which the coupled system is multi-stable. This multistability corresponds to the coexistence of mixed mode oscillations, bursting oscillations and subthreshold oscillations of the resonator neuron. In addition, some interval in the parameter space is introduced such that, for the values of the coupling strength within which the resonator neuron is in tonic spiking mode, while for the values of the coupling strength outside which the resonator neuron exhibits subthreshold oscillations. It is also verified that the final synchronization of the coupled neurons actually corresponds to the synchronization of tonic spiking oscillations of the integrator neuron and one-bursting oscillations of the resonator neuron.

\end{abstract}



\maketitle

\tableofcontents

\section{Introduction}
 Excitability is one of the most important characteristics of the neuron.
In 1948, Hodgkin identified three different types of excitability of the neuron. These types are qualitatively distinguished by the frequency-current relation. In the type
$I$ excitability the spiking frequency may be arbitrarily low depending on
the strength of the applied current, while in the type $II$ excitability the spiking
frequency is in a certain positive band \cite{hodgkin1948local}. 
Excitability and all of its types can be observed in the Hodgkin-Huxley model of neuronal dynamics (H-H model) as a four-dimensional system of differential equations \cite{hodgkin1952quantitative} and in at least two-dimensional reductions of H-H type models \cite{izhikevich2008dynamical}. Using the two-dimensional reductions of H-H type models is an efficient approach to studying the dynamics of the neuron. 

All the functions of the brain are obtained by  coupling of the neurons. Coupling the neurons is through electrical synapse (gap junction), chemical synapse or both of them. In the ninety decade, numerous studies on the effects of Connexins (Cx), i.e. gap junction proteins, on the genetic diseases and on the synchronization of the network of the neurons showed that the electrical synapse is one of the most important ways through which the signals of the neurons are transmitted \cite{hormuzdi2004electrical}. 
Investigating the dynamics of two coupled neurons through the gap junction may reveal some properties of networks of the neurons in brain.

Dynamics of two coupled oscillators have been investigated in science, especially in  natural sciences. In two coupled identical oscillators, in-phase oscillations, anti-phase oscillations and out-of-phase oscillations may be observed.
In addition, these types of oscillations may be observed in two coupled identical oscillatory neurons, which are coupled through the gap junction \cite{mirollo1990synchronization}.  Coupling two non-identical oscillatory neurons of the same types of excitability, may result in more complex dynamics such as, the existence of cascades of period-doubling bifurcations and coexistence of several stable limit cycles \cite{bindschadler2001bifurcation}. Investigating the dynamics of two weakly coupled identical neurons with the heterogeneity in synaptic conductances shows that, the solutions which correspond to in-phase oscillations and anti-phase oscillations coexist, but just one of them is stable \cite{baesens2013interaction}.

 More complex oscillation patterns may be observed in the  neurons.
Bursting is a dynamic state where a neuron repeatedly fires discrete groups or bursts of spikes. Each such burst is followed by a period of quiescence before the next burst occurs,  where periods of rapid action potential spiking are followed by quiescent periods much longer than typical inter-spike interval. Based on the dynamical systems theory, this type of oscillations may be observed in at least three-dimensional models \cite{coombes2005bursting,izhikevich2000neural,rinzel1987formal}. 
Burst synchronization of the coupled  neurons  is typically used to refer to a temporal relationship between active phase onset or offset times across the neurons. That is, the neurons  start the firing and become silent almost simultaneously. Burst synchronization may be observed in the brain \cite{kamioka1996spontaneous,steriade1993thalamocortical}. Burst synchronization may be also observed in two coupled neurons of different types of  excitability \cite{de2020burst}. The results of \cite{razvan2020} show that  two coupled integrator neurons of different types of  excitability may exhibit  burst synchronization, when none of the neurons burst intrinsically. 
 
 Another oscillation pattern  which is reported in  nature \cite{petrov1992mixed,zhabotinskii1964periodic} and  in the neuron \cite{drover2004analysis,epstein1996nonlinear,koper1995bifurcations,medvedev2004multimodal,rotstein2006reduced} is mixed mode oscillations. Mixed mode oscillations (MMOs) is an oscillation pattern within which the system exhibits a combination of small and large amplitude oscillations.
By using the results of \cite{krupa2001relaxation} one can conclude that   mixed mode oscillations may be observed in at least three-dimensional autonomous ordinary differential equations.

A great deal of research has been devoted on the investigation of neuronal dynamics in coupled neurons.  One of  the  most important questions is that, in a network of  neurons, how the differences  in the   dynamics of the neurons  affects the dynamics of the network. In \cite{razvan2020} the dynamics of two coupled integrator neurons of different types of excitability through the gap junction has been investigated. The results of \cite{razvan2020} shows that, in a network of  coupled neurons of different types of excitability the neurons may exhibit oscillations against their types of excitability.

Some types of neurons have the tendency to fire at particular frequencies, so-called resonators \cite{izhikevich2008dynamical}.  In the resonator neurons oscillatory activity can  be observed in the form of subthreshold membrane potential oscillations (i.e. in the absence of action potentials) \cite{llinas1986oscillatory}. 
  In a network of neurons, the dynamics of resonator neurons is more related to the frequency of input signal to the resonator neurons rather than  the coupling strength.
  The results of \cite{muresan2007resonance} show that resonance and integration at the neuron level might interact in the brain to promote stability as well as flexibility and responsiveness to external input and that membrane properties, in general, are essential for determining the behaviour of large networks of neurons. In \cite{schwabedal2016qualitative}  all oscillation patterns of the network of  three coupled resonator neurons have been investigated, where the individual neurons in the network are burster neurons.   
  
   Another important question is that how differences in intrinsic characteristics of the  coupled neuron's dynamics, e.g. resonance or integratory dynamics,  affects the network's dynamics.  In this manuscript, the answer to this question is investigated  when a silent resonator neuron is coupled with a spiking integrator neuron through the  gap junction. More precisely, $I_{Na,K}$ model  \cite{izhikevich2008dynamical}, which is equivalent to the famous and widely used Morris-Lecar $I_{Ca}+I_K$ model \cite{morris1981voltage}, is considered as a model which describes the neuron's dynamics. Hence, none of  the coupled neurons exhibit mixed mode oscillations and bursting oscillations intrinsically.  Then, by using the dynamical systems theory (e.g. the bifurcation theory) it is examined that, how increasing the coupling strength affects the dynamics of the neurons, when one of the coupling strength is fixed and the other varies. In this manuscript all the observed oscillation patterns of the neurons are introduced, then the transition mechanisms between different patterns are examined.

More specifically, the coupled system is initially bistable, a stable limit cycle which corresponds to the subthreshold oscillations  of the resonator neuron and a stable two-dimensional torus which corresponds to the phase-locking oscillations of the neurons. By increasing the coupling strength, different oscillation patterns may be observed in the resonator neuron, while the integrator neuron is in tonic spiking mode.
The observed patterns in the resonator neuron are, phase-locking oscillations of the neurons, mixed mode oscillations (MMOs), bursting oscillations, subthreshold oscillations, intermittent oscillations which is a combination of spiking oscillations and subthreshold oscillations and synchronous oscillations of the neurons. 
By increasing the coupling strength,  multi-stability, as one of the most important features of the coupled system, is observed  in the system. That is, for strong enough coupling strength the coupled system  has three stable limit cycles,  which correspond to the mixed mode oscillations,  bursting oscillations   and   subthreshold oscillations of the resonator neuron. 
The bursting oscillations which is  observed in the resonator neuron,  have  a feature  that distinguishes them from the other  observed bursting oscillations. Usually, in    the  bursting oscillations the stable limit  cycle, which corresponds to the subthreshold oscillations, loses its stability.  It is noticeable that in our coupled system the stability of the limit cycle which corresponds to the subthreshold oscillations persists.   It is  also stated that 
 in  the coupled system   the final synchronization of the neurons actually corresponds to the synchronization of tonic spiking oscillations of the integrator neuron and 1-bursting oscillations of the resonator neuron.
 The achievements  of this manuscript also confirm that   the resonator neurons may fire through the fine tuning. More precisely, some  interval in the parameter space is introduced  such that for the  values of the coupling strength  within this interval the resonator neuron is in spiking mode, while for the values of the coupling strength  outside which   the resonator neuron exhibits the subthreshold oscillations. 

The  manuscript is organized as follows: 

In   subsection $2.1$, the $I_{Na,K}$ model    is concerned as a model which describes the neuron's dynamics. Then, the parameters of the model is chosen such that the model exhibits types $I$ and $II$ excitability. In   section $2.2$, the coupled system  is introduced. In   section $3$, all the observed oscillation patterns in the resonator neuron are introduced, while the integrator neuron is in tonic spiking, more specifically,  phase-locking oscillations of the neurons in the subsection $3.1$, emergent mixed mode oscillations (MMOs)  of the  resonator neuron in the subsection $3.2$, bursting oscillations of the resonator neuron in the subsection $3.3$, subthreshold oscillations of the  resonator neuron  in the subsection $3.4$, intermittent oscillations of the  resonator neuron in the subsection $3.5$, synchronous oscillations of the   neurons in the subsection $3.6$. Then,    by using the dynamical systems theory (e.g. the bifurcation theory), it is investigated that how increasing the coupling strength affects each oscillation pattern.
In  section $4$, transition mechanisms between different  oscillation patterns are examined. More precisely,  transition from phase-locking oscillations to MMOs in the subsection $4.1$,  transition from MMOs to  the bursting oscillations in the subsection $4.2$, transition from bursting oscillations to the subthreshold oscillations in the subsection $4.3$, transition from subthreshold oscillations to intermittent oscillations in the subsection $4.4$ and the transition from the intermittent oscillations to synchronous oscillations in the subsection $4.5$, are investigated.

Note that all simulations in this manuscript are done by XPPAUT \cite{ermentrout2002simulating}.

\section{The model}

Excitability is one of the most important characteristics of a neuron. In 1948,  Hodgkin identified   three different  types of excitability of the neurons. These types are qualitatively distinguished by the frequency-current relation. In type
$I$ excitability, the spiking frequency may be arbitrarily low depending on
the strength of the applied current, while in type $II$ excitability the spiking
frequency is in a certain positive band \cite{hodgkin1948local}. 
Excitability and all of its types can be observed in Hodgkin-Huxley model  of neuronal dynamics (H-H model) as a  four-dimensional system of differential equations \cite{hodgkin1952quantitative} and in  at least   two-dimensional reductions of  H-H type models \cite{izhikevich2008dynamical}. In the  following, it is stated that $I_{Na,K}$ model, which is equivalent to the famous and widely used Morris-Lecar $I_{Ca}+I_K$ model \cite{morris1981voltage},
 shows two different types of excitability for different values of the parameters. Then, the  system of two coupled neurons of different types of excitability (the coupled system) is introduced.

\subsection{\textbf{The $I_{Na,K}$ model and different types of excitability}}
In this manuscript,   $I_{Na,K}$ model is considered as a model of single neuron's dynamics \cite{izhikevich2008dynamical}. $I_{Na,K}$ model, which is equivalent to the famous and widely used Morris-Lecar $I_{Ca}+I_K$ model  \cite{ morris1981voltage}, 
 consists of a fast
$Na^+$ current and a relatively slower $K^+$ current. A reasonable assumption based on experimental observations is that, the $Na^+$ gating
variable $m(t)$ is much faster than the voltage variable $V (t)$, so that $m$ approaches the
asymptotic value $m_{\infty} (V )$ practically instantaneously. In this case, one can substitute
$m = m_{\infty} (V )$ into the voltage equation \cite{izhikevich2008dynamical}.   The $I_{Na,K}$ model is as follows:
\begin{equation}\label{1}
\left\{\begin{array}{l}
\dot{V}=\dfrac{-(g_L(V-E_L)-g_{Na}m_{\infty}(V)(V-E_{Na})-g_Kn(V-E_K))}{C}+I, \\
 \dot{n}=\dfrac{(n_{\infty(V)}-n)}{\tau(V)},
\end{array}\right.
\end{equation}

\begin{figure} [ht]
\begin{center} 
\includegraphics[width=.5\linewidth,height=.07\linewidth]{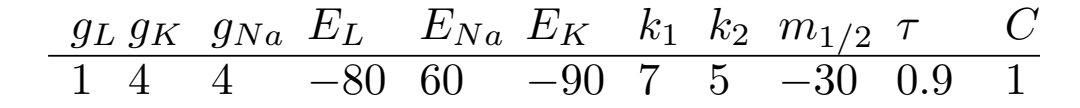}
 \label{firsttable}
       \caption{The fixed parameters values of the system  (\eqref{1}).}
\end{center}
\end{figure}

where
 \begin{center}
$m_{\infty}(V)=\Big(1+\exp(\frac{m_{1/2}-V}{k_1})\Big)^{-1},$
$n_{\infty}(V)=\Big(1+\exp(\frac{m'_{1/2}-V}{k_2})\Big)^{-1},$ 
\end{center} 
and all the fixed parameters are listed in Figure 1.
In the following, it is shown that the system (\eqref{1}) displays different types of excitability for different values of $m'_{1/2}$.

Consider a system of ordinary differential equations $\dot{x}=f(x), \hspace{2mm} x\in R^n$ with an equilibrium $x_0$ (that is, $f(x_0)=0$). $x_0$ is called a hyperbolic point if  all the eigenvalues of Jacobian matrix at $x_0$ have nonzero real part.   $x_0$ is stable if   all the eigenvalues of Jacobian matrix at $x_0$ have negative real parts. Moreover, $x_0$  is called  a saddle point if some of the eigenvalues of  Jacobian matrix have positive real parts and the others have negative real parts \cite{kuznetsov2013elements}.

 Now, suppose that the system has a limit cycle $L_0$  passing through $x_0$. 
By choosing a suitable Poincar\'e section to the limit cycle $L_0$, one can define the corresponding Poincar\'e  map, $T$.   The fixed point of the Poincar\'e map $T$ is obtained by numerically solving 
$T(B)-B=0$.
The stability of the fixed point is governed by the eigenvalues of the Jacobian matrix $[\partial T/\partial B]$. A local bifurcation of the periodic solution occurs when the Jacobian matrix evaluated at the fixed point has an eigenvalue of the absolute value of one, $\mu$. 
Suppose that $\mu$ is the eigenvalue characterises the type of local bifurcation \cite{kuznetsov2013elements}. When $\mu=1$ then the fixed point of the Poincar\'e map corresponds to a saddle-node (SN) bifurcation point, when $\mu=-1$ it corresponds to a period-doubling bifurcation point and  when $\mu=e^{i\theta}$ with $\theta \neq 0$,  it corresponds to  a Neimark-Sacker bifurcation point.

\begin{remark}
 The stable equilibrium of the system  (\eqref{1}) corresponds to the silent phase of the neuron. In addition, the stable  limit cycle of   (\eqref{1}) corresponds to the spiking oscillations of the neuron.
\end{remark}

\begin{figure} [!ht]
\begin{center} 
\includegraphics[width=1\linewidth,height=.95\linewidth]{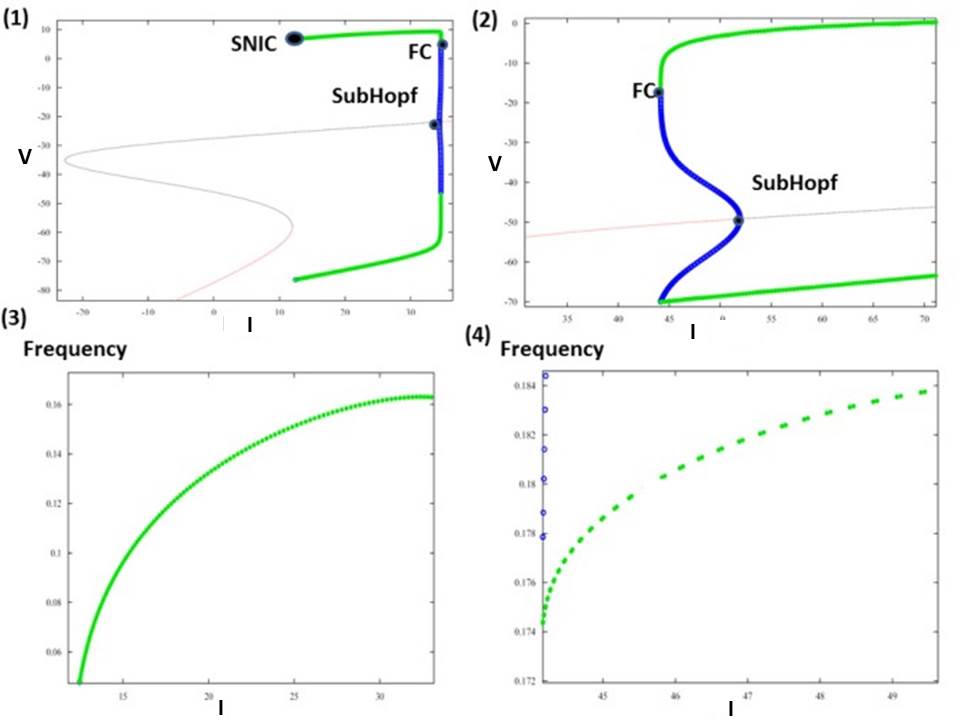}
\caption{$(1),(2)$. Bifurcation diagram of the system (\eqref{1}) for two different values
$m'_{1/2}$
when
$I$ is the bifurcation parameter.  Here, the color red indicates the stable equilibrium,
the color green indicates the stable limit cycle  and the blue color indicates the unstable limit cycle of the system (\eqref{1}). (1) For $m'_{1/2}=-30$, by increasing I,
stable equilibrium of the system  undergoes Saddle-Node bifurcation on Invariant Circle ($SNIC$ bifurcation), then a stable limit cycle, $W_1$, appears. (2) For $m'_{1/2}=-45$, stable equilibrium of the system
undergoes Subcritical Hopf bifurcation,  hence it becomes unstable. Then the state of the system tends to the  stable limit cycle $W_2$.  (3),(4).  The frequency-current relation. (3) For
$m'_{1/2}=-30$, spiking frequency can be arbitrarily low depending on the strength of the applied
current, i.e. the system shows type $I$ excitability. (4) For $m'_{1/2}=-45$, the spiking frequency
is in a certain positive band depending on the strength of the applied current, i.e. the
system shows type $II$ excitability. The fixed parameters values of the system (\eqref{1}) have been listed in the
Figure . } 
\label{bif-fri}
\end{center}
\end{figure}
Let $m'_{1/2}=-30$. Then the system (\eqref{1}) has a stable equilibrium, a saddle equilibrium and an unstable one (Figure \ref{bif-fri} (1)). The stable equilibrium corresponds to the silent phase of the neuron. By increasing $I$, the stable equilibrium and the saddle one get closer to each other. By further increasing in $I$, the stable equilibrium and the saddle one coalesce at $I=11.99$.  that is, the heteroclinic orbit between them becomes a homoclinic orbit. By increasing $I$ the saddle-node equilibrium disappears and a stable limit cycle $W_1$, which surrounds the unstable equilibrium, appears. In the other words, the system undergoes saddle-node bifurcation on invariant circle ($SNIC$ bifurcation). 
Through this bifurcation the neuron initiates periodic spiking.

\begin{remark}
After the $SNIC$ bifurcation, 
 tonic spiking starts from the zero frequency (Figure \ref{bif-fri} (3)), That is, for  $m'_{1/2}=-30$ the system shows type $I$ excitability.
\end{remark}
Now, let $m'_{1/2}=-45$. Then the system has two stable invariant sets, a stable equilibrium $A$ which corresponds to the silent phase of the neuron and a stable limit sycle $W_2$ which corresponds to tonic spiking mode in the neuron (Figure \ref{bif-fri} (2)). The basins of attraction of these invariant sets are separated by an unstable limit cycle.  As $I$ increases, the  unstable limit cycle $W_3$ gets smaller,  shrinks to  $A$ and finally makes it lose stability. 
 That is, at $I=51.9$ the system undergoes subcritical Hopf bifurcation ("SubH" bifurcation). As a result, the state of the system tends to the stable limit cycle $W_2$ and the neuron initiates spiking oscillations. 
\begin{remark}
After  "SubH" bifurcation, the neuron initiates the spiking oscillations with a positive frequency. In the other words, for $m'_{1/2}=-45$ the system shows type $II$ excitability.    (Figure \ref{bif-fri} (4))
\end{remark}

\subsection{\textbf{The Model}}

\begin{figure}[!ht]
\begin{center}
\includegraphics[width=1\linewidth,height=.7\linewidth]{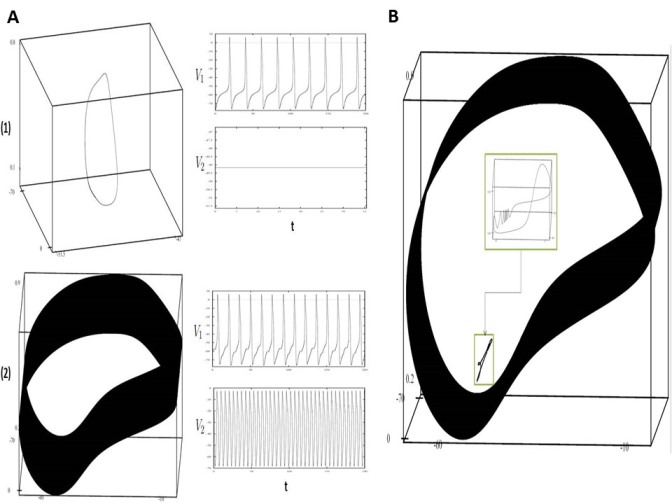} 
\caption[] {When $q_1 = q_2 = 0$  the coupled system \eqref{2} is bistable, a stable limit cycle  and a
stable two-dimensional torus.  (1),(2).  Three-dimensional  image, $(V_2,V_1,n_2)$,   and the corresponding voltage time series of $"I"$ and $"II"$ for $(1)$ the stable limit cycle  and   $(2)$ the stable two-dimensional torus.}
\label{Figure2}
\end{center} 
\end{figure}

\begin{figure} [ht]
\begin{center} 
\includegraphics[width=.7\linewidth,height=.08\linewidth]{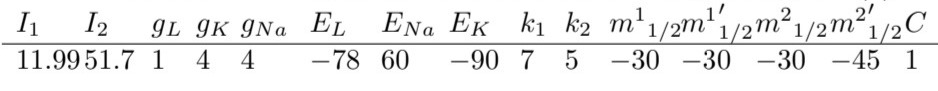}
 \label{fixedparamstwo}
       \caption{The fixed parameters values of the system  \eqref{2}.}
\end{center}
\end{figure}

In this article, a spiking neuron of type $I$ excitability (neuron $"I"$, $m'_{1/2}=-30$) is coupled with a silent neuron of type $II$ excitability (neuron $"II"$, $m'_{1/2}=-45$). Coupling between the neurons is through a linear form of gap
junction, $q_j(V_i -V_j)$, where $V_i$ and $V_j$ are the voltage of the neurons,  $q_j$
is the coupling strength and $i,j \in\{1,2\}$. The coupled system is as follows:
 \begin{equation}\label{2}
\left\{\begin{array}{l}
 \dot{V_1}=I_{total}(V_1,n_1)+I_1+q_1(V_2-V_1), \\
\dot{n_1}=\frac{n_{\infty}(V_1)-n_1}{\tau_1},  \\
 \dot{V_2}=I_{total}(V_2,n_2)+I_2+q_2(V_1-V_2), \\
 \dot{n_2}=\frac{n_{\infty}(V_2)-n_2}{\tau_2},\\
\end{array}\right.
\end{equation}

where

  $$I_{total}(V,n)=\frac{-(g_L(V-E_L)-g_{Na}m_{\infty}(V)(V-E_{Na})-g_Kn(V-E_K))}{C},$$
and the fixed values of the parameters of (\eqref{2}) are listed in   Figure 4.
 
\begin{figure}[ht]
\begin{center}
\includegraphics[width=1\linewidth]{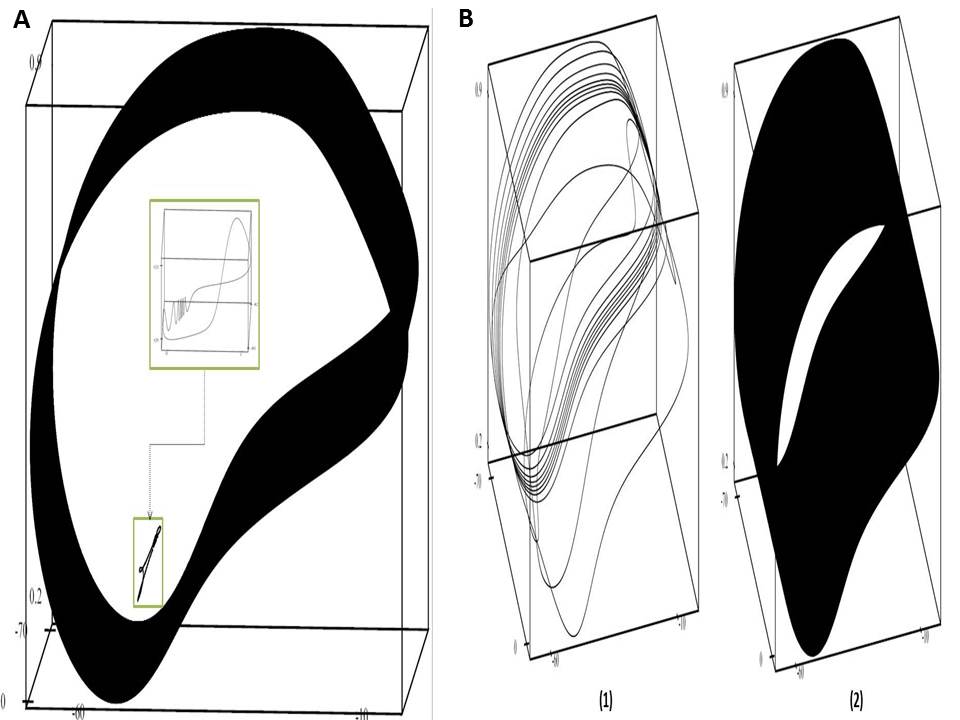} 
\caption[]{A. Three-dimensional image, $(V_1,V_2,n_2)$,  of the bistability, i.e. co-existance of  a stable limit cycle and a  stable two-dimensional torus for $q_{1} =0.05$ and $q_{2}=0.04$. B. The dynamics on the smooth torus. (1) periodic dynamics for $q_2=0.082$ and  (2) quasi-periodic dynamics for $q_2=0.03$. }
\label{Figure3}
\end{center} 
\end{figure}

When $q_1 = q_2 = 0$,  the  system (\eqref{2}) is bistable, a stable limit cycle $W_1\times \{A\}$ and a
stable two-dimensional torus $W_1\times W_2$. Moreover, the system has a saddle two-dimensional torus $W_1\times W_3$, whose stable manifold separates the basins of attraction of
the stable limit cycle and the stable torus (Figure \ref{Figure2}).

For  $q_1 = q_2 = 0$,
all the invariant sets of the coupled system are hyperbolic. Hence,   the
conditions of structural stability hold for the invariant subsets of the system (\eqref{2})
\cite{hirsch1970invariant}. As an immediate consequence of structural stability, bistability of the
system persists for small values of $q_1$ and $q_2$. More precisely, for small values of $q_1$ and $q_2$  the coupled system is bistable, a stable limit cycle which corresponds to tonic spiking oscillations of  $"I"$ and subthreshold oscillations of  $"II"$, and a stable two-dimensional torus which corresponds to phase-locking
 oscillations  of $"I"$ and $"II"$ (Figure \ref{Figure3} A). Moreover, the system has a saddle two-dimensional torus  whose stable manifold separates the basins of attraction of
the stable limit cycle and the stable torus.      Throughout this
manuscript, we fix $q_1 = 0.05$ and then we investigate the dynamics of the coupled system $(\eqref{2})$, when $q_2$ varies. Here, all the fixed parameter's values of the coupled system (\eqref{2}) are listed in Figure 4.

 \section{\textbf{Different observed  oscillation patterns in the coupled neurons}}

 As the coupling strength increases, different  oscillation patterns are observed in the coupled  neurons.  More specifically, the observed oscillation patterns in one of the coupled neurons or both of them are as follows: Phase-locking oscillations of the neurons,  mixed mode oscillations of the resonator neuron,  bursting oscillations of the resonator neuron, subthreshold oscillations of the resonator neuron, intermittent oscillations of the resonator neuron and synchronous oscillations of the neurons. Notice that, none  of the coupled  neurons $"I"$ and $"II"$, as two-dimensional ordinary differential equations,  can  intrinsically show some of these patterns, such as  mixed mode oscillations and  bursting oscillations. In the following, all the observed  patterns are introduced, then the effect of increasing the coupling strength on each of these patterns is examined.
 \subsection{\textbf{Phase-locking oscillations}}
 As mentioned in the section $1.2$, for sufficiently small values of $q_2$ the coupled system (\eqref{2}) is  bistable. That is, the system has two stable invariant sets, a stable limit cycle which corresponds to tonic spiking oscillations of  $"I"$ and subthreshold oscillations of  $"II"$, and a stable smooth two-dimensional torus which corresponds to phase-locking
 oscillations of  $"I"$ and $"II"$  (Figure \ref{Figure3} A).

 Based on the Poincar\'e-Denjoy theory,  the dynamics on the smooth two-dimensional  torus is determined by the  rotation number which corresponds to a suitable Poincar\'e map.  More precisely,  when the torus is smooth,  one may choose  a suitable Poincar\'e section to the torus such that, the corresponding  Poincar\'e map  has a domain  which is topologically a circle, $S^1$. Hence,  the rotation number is defined for the Poincar\'e map. When the  rotation number is rational, the dynamics on the smooth torus is periodic \cite{kuznetsov2013elements}. Moreover, the dynamics on the torus is quasi-periodic, when the rotation number is irrational \cite{denjoy1932courbes}.

For suitable values of the coupling strengths,  phase-locking oscillations 
 is observed in  $"I"$ and $"II"$.  Figure  \ref{Figure3} B shows  three-dimensional images of periodic and quasi-periodic dynamics on the torus. 

 \subsection{\textbf{Tonic spiking / emergent mixed mode oscillations (TS/MMOs)}} 
 In dynamical systems, mixed mode oscillations (MMOs) is an   oscillation pattern  within which the system shows  oscillations with two or more distinct lengths \cite{albahadily1989mixed,del2002periodicity,zhabotinskii1964periodic}. Based on the theory of dynamical systems, this type of oscillations may be observed in at least three-dimensional models \cite{koper1995bifurcations,milik1998geometry}.  Hence, none of the neurons $"I"$ and $"II"$, as two-dimensional ordinary differential equations,  can   exhibit mixed mode oscillations intrinsically.

\begin{figure}[!ht]
\begin{center}
\includegraphics[width=1\linewidth]{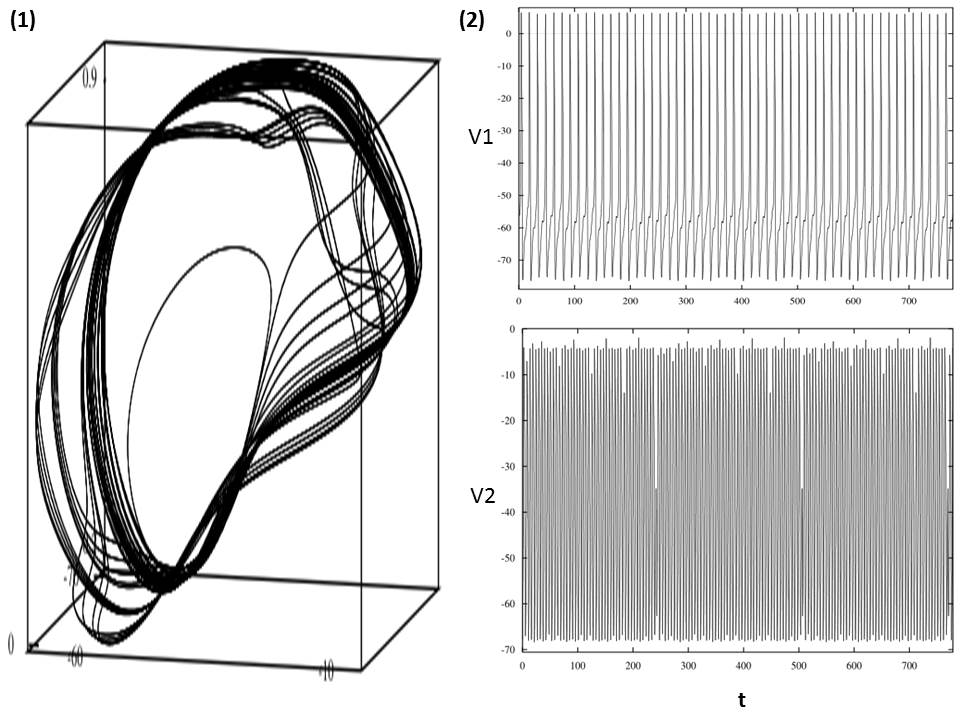} 
\caption[]{(1) Three-dimensional image, $(V_1,V_2,n_2)$, of  the stable limit cycle $M$, which corresponds to the mixed mode oscillations of $"II"$, for $q_{2}=0.086783$.   (2)  voltage time series of $"I"$ and $"II"$,  $V_1$ and $V_2$ respectively, corresponding to $M$.}
\label{Figure4}
\end{center} 
\end{figure}

By increasing the coupling strength, mixed mode oscillations is observed in $"II"$  while $"I"$  is in tonic spiking  mode.  More precisely, by increasing the coupling strength, for suitable initial conditions, the state of the system approaches to a stable limit cycle $M$ (Figure  \ref{Figure4} (1)).   As shown in Figure  \ref{Figure4} (2), the voltage time series of "II" shows two distinct oscillations with different lengths. That is, $M$ corresponds to the mixed mode oscillations of "II", while "I" is in the tonic spiking mode. 
In the following, it is investigated that how increasing $q_{2}$ affects mixed mode oscillations of "II". 

\begin{figure}[h!]
\begin{center}
\includegraphics[width=1\linewidth]{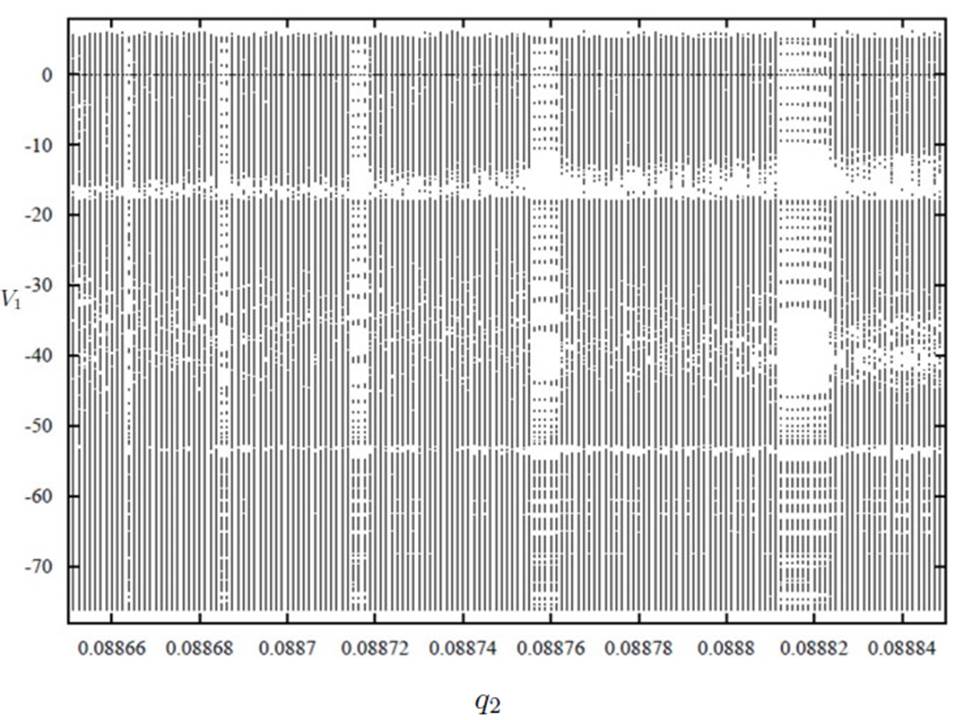} 
\caption[]{Bifurcation diagram of the Poincar\'e map of the coupled  system (\eqref{2}) corresponding to the section $V_2=-50$, where $q_2$ is the bifurcation parameter. By increasing the coupling strength, the system (\eqref{2}) undergoes cascades of period-doubling bifurcations.}
\label{Figure5}
\end{center} 
\end{figure} 

\begin{figure}[h!]
\begin{center}
\includegraphics[width=1\linewidth]{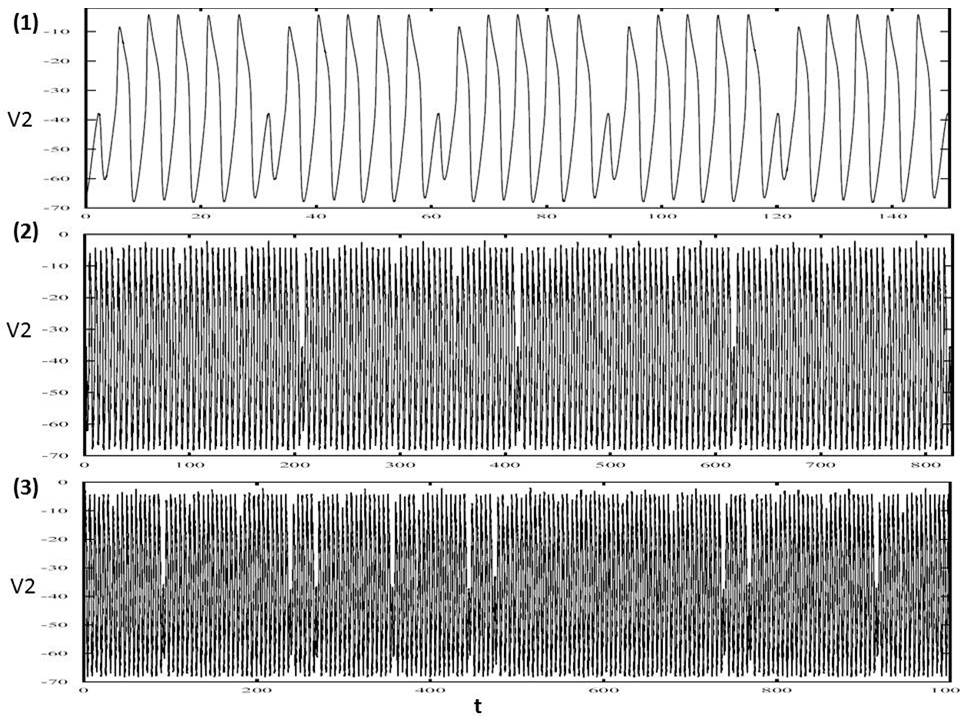} 
\caption[]{(1) Simple and periodic MMOs for $q_2=0.092$ of  $"II"$, (2) Complex and  periodic MMOs for  $q_2=0.08872$ of  $"II"$, (3) Chaotic MMOs for  $q_2=0.093$ of  $"II"$.}
\label{Figure6}
\end{center} 
\end{figure} 

Sander and York in \cite{sander2012connecting} stated that:

\begin{theorem}
Assume that the following conditions hold for the smooth map
$F(x,\mu)$, where $x \in  R^n$ and $\mu \in R$. Suppose that $\mu_1$ and $\mu_2$ are two parameter
values and $\mu_1<\mu_2$:

1. $F$ is infinitely differentiable in $\mu$ and $x$, and all of its bifurcations
(including saddle-node, Hopf and period-doubling bifurcation) are generic.

2. All periodic orbits at $\mu_1$ and $\mu_2$ are hyperbolic.

3. The number of periodic orbits at $\mu_1$ is finite.

4.There is a number $G>1$  at $\mu_2$ for which the number of periodic orbits
of period $P$ at $\mu_2$, $fixed(\mu_2; P)$, satisfies: $fixed(\mu_2; P) \geq G^P$ for infinitely
many $P$.

5. All but a finite number of periodic orbits at $\mu= \mu_2$ have the same
unstable dimension.

Then, there are infinitely many distinct period-doubling cascades between
$\mu_1$ and $\mu_2$. In addition, in each of these cascades, the chaotic regime ends
with homoclinic bifurcation to fold limit cycle.
\end{theorem}

Figure \ref{Figure5} shows the bifurcation diagram of the Poincar\'e map of the system (\eqref{2}) in $[0.08868 , 0.08884]$.
As demonstrated by this   diagram, increasing the coupling strength makes  the system undergo cascades of period-doubling bifurcations.  Now, by using the Theorem 3.1,  one can conclude that: 
 
 \begin{proposition}
 By increasing $q_2$ in the interval  $[0.08868 , 0.08884]$, the coupled system undergoes infinitely many  cascades of period-doubling bifurcations. Each of these cascades ends up through homoclinic bifurcation to a saddle-node cycle.
 \end{proposition}

\begin{proof}
The assumptions 1--5 of the theorem 3.1 are clearly established. As demonstrated by the bifurcation diagram of the Poincar\'e map of the system (\eqref{2}) (Figure \ref{Figure5}), by increasing $q_{2}$ the limit cycle M undergoes cascades of period-doubling bifurcations. For instance, at $q_{2}= 0.0886856$ the cycle has a floquet multiplier $\lambda = -1$, then the system undergoes 
period-doubling bifurcation. The cascade eventually leads to the chaotic dynamics. More precisely, at $q_{2}= 0.08869$ the dynamics is chaotic. Since the system undergoes cascades of period-doubling bifurcation, the number of the  limit cycles of the system grows exponentially. Hence, the assumption 6 of the theorem 1 also holds.

Now, by using the theorem 3.1 one can  conclude that, there are infinitely many distinct   period-doubling cascades  in $[0.08868 , 0.08884]$.
Each of these cascades leads to   chaotic dynamics. Moreover, in each cascade, the chaotic oscillations ends up through homoclinic bifurcation to saddle-node cycle, then a stable limit cycle appears.
\end{proof}

 By further increasing in $q_2$, simple periodic MMOs (Figure \ref{Figure6} (1)), Complex periodic MMOs (Figure \ref{Figure6} (2)), and chaotic MMOs (Figure \ref{Figure6} (3)) are observed in $"II"$, while $"I"$ is in tonic spiking mode.

\begin{figure}[!ht]
\begin{center}
\includegraphics[width=1\linewidth]{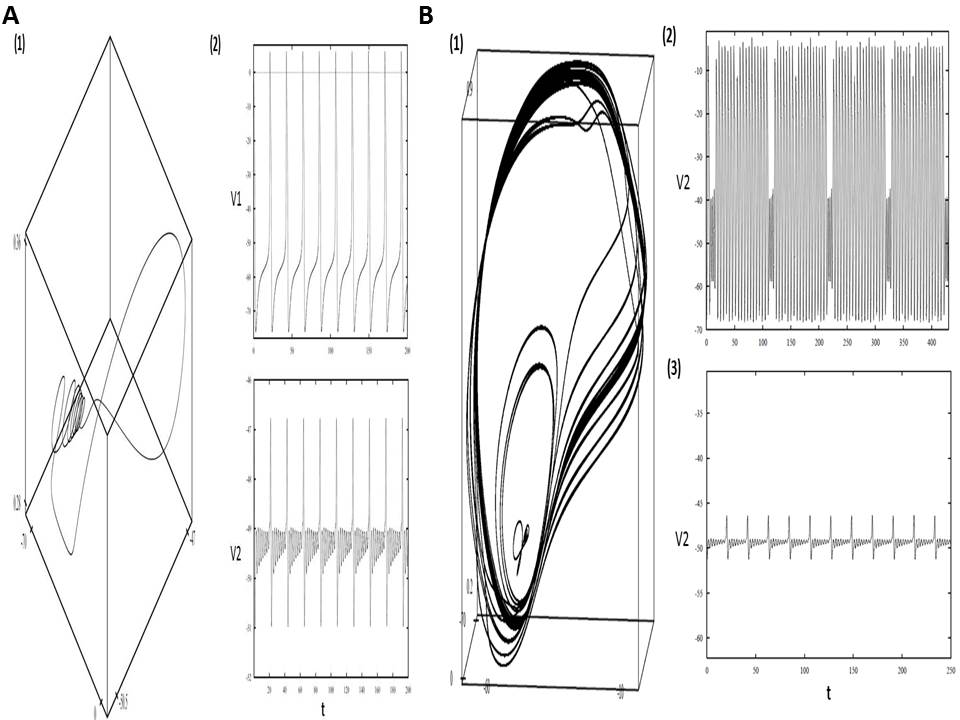} 
\caption[]{A.   By increasing the coupling strength, the stable limit cycle which corresponds to the subthreshold oscillations of $"II"$ persists. (1) Three-dimensional image, $(V_1,V_2,n_2)$, of the stable limit cycle and (2) the  corresponding voltage time series of  $"I"$ and $"II"$ for  $q_2=0.08675$. B. For $q_2=0.099097$, the coupled system has two stable limit cycles. One of which corresponds to the bursting oscillations of $"II"$ and the other corresponds to the subthreshold oscillations of $"II"$. (1) Three-dimensional image of the bistability. The voltage time series of the resonator neuron which corresponds to  (2) the bursting oscillations and (3) the subthreshold oscillations.  }
\label{Figure7}
\end{center} 
\end{figure}

\begin{remark}
It is noticeable that, for suitable initial conditions subthreshold oscillations is observed in $"II"$ while $"I"$ is in tonic spiking mode (Figure \ref{Figure7} A). More precisely, after the torus destruction the coupled system has two stable limit cycles. One of them corresponds to the subthreshold oscillations of $"II"$, and the other corresponds to the mixed mode oscillations of $"II"$. Moreover, the coupled system (\eqref{2}) has an unstable torus whose stable manifold separates the basins of attraction of the stable limit cycles.
\end{remark}
\begin{remark}
 In the coupled system (\eqref{2}), mixed mode oscillations are actually canard-induced MMOs \cite{krupa2001relaxation}.  Let $\Gamma$ be the trajectory of the system which corresponds  to the MMOs. $\Gamma$ sometimes moves along the stable manifold of the unstable torus toward  the unstable torus. Hence,     small amplitude oscillations is observed in $"II"$. After a while, $\Gamma$ returns to a neighbourhood of the stable limit cycle, hence   spiking oscillations of $"II"$ is observed again. 
\end{remark}

 \subsection{\textbf{Tonic Spiking/Bursting Oscillations (TS/Bursting)}}
 Bursting is a dynamic state from slow oscillations from a low voltage (silent phase) to a plateau upon which is superimposed rapid spiking (active phase) \cite{izhikevich2000neural}. This type of oscillations may be observed in at least three-dimensional models. Hence, none of the coupled neurons $"I"$ and $"II"$, as a two-dimensional model, are able to exhibit bursting oscillations intrinsically.

   \begin{figure}[!ht]
\begin{center}
\includegraphics[width=1\linewidth]{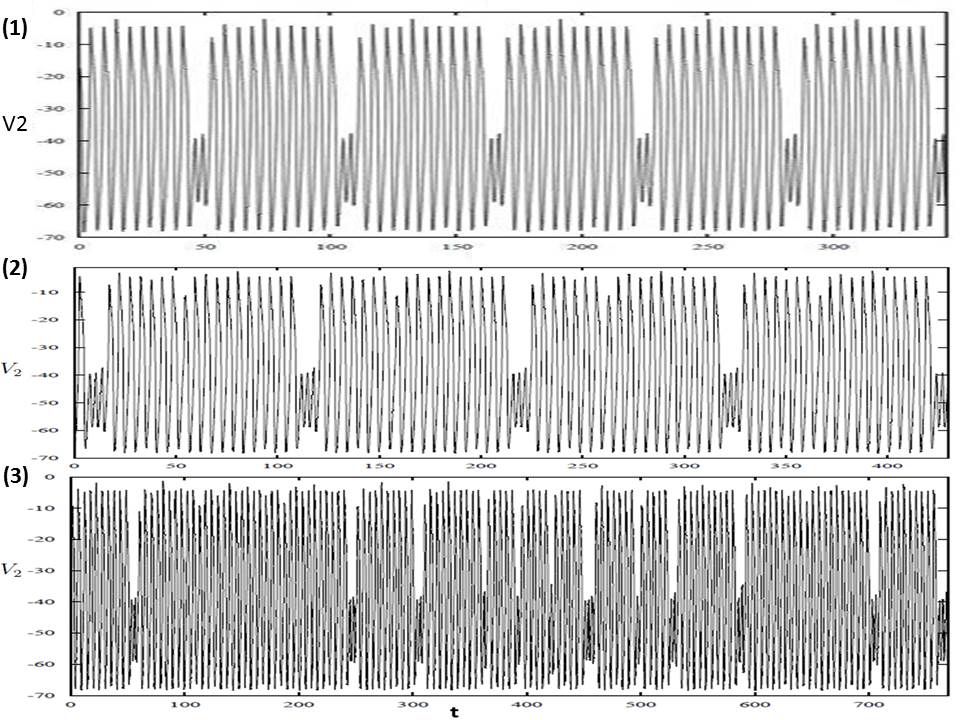} 
\caption[]{ Bursting oscillations of
$"II"$.  Periodic bursting oscillations for (1) $q_2=0.096783$ and (2) $q_2=0.099097$. (3) chaotic bursting oscillations for
$q_2=0.09919$.}
\label{Figure8}
\end{center} 
\end{figure} 

 \begin{figure}[!ht]
\begin{center}
\includegraphics[width=1\linewidth]{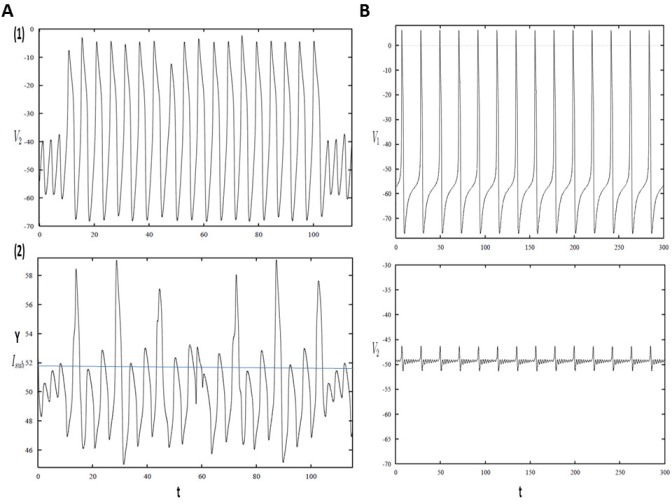} 
\caption[]{ A. For $q_2=0.099097$, bursting oscillations is observed in $"II"$. (1) The voltage
time series of $"II"$, $V_2$, corresponding to the bursting oscillations.
(2) The corresponding input signal $Y$ to $"II"$ through the input
current and the coupling strength, i.e. $Y = I_2 +q_2(V_1-V_2)$, during the active phase of $"II"$. B. For strong enough coupling strength, the resonator neuron exhibits subthreshold oscillations. The corresponding voltage
time series of (1) $"I"$ and (2) $"II"$ for $q_2=0.1019$.}
\label{Figure9}
\end{center} 
\end{figure}

As  previously mentioned, after the torus destruction, the coupled system has two stable limit cycles, $U_1$ and $U_2$. The stable limit cycle  $U_1$ corresponds to the subthreshold oscillations of $"II"$, and $U_2$ corresponds to the mixed mode oscillations of $"II"$. Moreover, the system has an unstable torus. The stable manifold of the unstable torus separates the basins of attraction of the stable limit cycles. 
 By further increasing in the coupling strength,  the unstable torus   gets more smaller.  Hence, starting from the points in the basin of attraction of  $U_2$, the trajectories spend more times in a neighbourhood of $U_1$. That is, $"II"$ exhibits  bursting oscillations. Some of the observed  bursting oscillations have been depicted in Figure \ref{Figure8}.

 \begin{remark}
 In the coupled system (\eqref{2}), the bursting oscillations  have  a feature  that distinguishes them from  most of the other  observed bursting oscillations. In   the bursting oscillations, usually   the stable limit  cycle, which corresponds to the subthreshold oscillations, loses its stability, though in the coupled system (\eqref{2}) the stability of the limit cycle persists. 
 
  \end{remark}
  In the following,  the observed bursting oscillations for $q_2=0.099097$ is investigated precisely:

Let $q_2=0.099097$. Then,  the coupled system (\eqref{2}) has two stable limit cycles (Figure \ref{Figure7} B(1)). One of them corresponds to the subthreshold oscillations of $"II"$ (Figure \ref{Figure7} B(3)). In the following, it is examined that the other stable limit cycle corresponds to the bursting oscillations of $"II"$ (Figure \ref{Figure7} B(2)). 
 
Let $I_{sub}:=51.9$. As mentioned in the section $2.1$, at $I=I_{sub}$, $"II"$ undergoes subcritical Hopf bifurcation, then $"II"$ initiates spiking oscillations.
Now, let $Y$ be the input signal to  $"II"$ through the input current and   the coupling  strength which is obtained by
 $Y=I_2+q_2(V_1-V_2).$
 As shown by Figure \ref{Figure9} (2), there exist some time intervals within which input signal to   $"II"$ is below the corresponding bifurcation value $I_{sub}$, hence in each of these time intervals  $"II"$ exhibits subthreshold oscillations. Going through each of these time intervals, the  input signal to  $"II"$  becomes bigger than the corresponding bifurcation value $I_{sub}$, therefore $"II"$ initiates periodic spiking  (Figure \ref{Figure9} (1)).  The periodic spiking continues until that at the next time interval $"II"$ becomes silent again. In the other words,  $"II"$ shows bursting oscillations. Since the input signal to $"II"$ is periodic, for $q_{2}=0.099097$ periodic bursting oscillations is observed in $"II"$.
 
 \subsection{\textbf{Tonic Spiking / Subthreshold Oscillations (TS/Sub)}}
 Some types of neurons have the tendency to fire at particular frequencies, so-called resonators \cite{izhikevich2008dynamical}.  In the resonator neurons, oscillatory activity can also be observed in the form of subthreshold membrane potential oscillations (i.e. in the absence of action potentials) \cite{llinas1986oscillatory}. 
  In a network of neurons, the dynamics of resonator neurons is more related to the frequency of input signal to the resonator neurons rather than  the coupling strength.

 \begin{figure}[!ht]
\begin{center}
\includegraphics[width=1\linewidth]{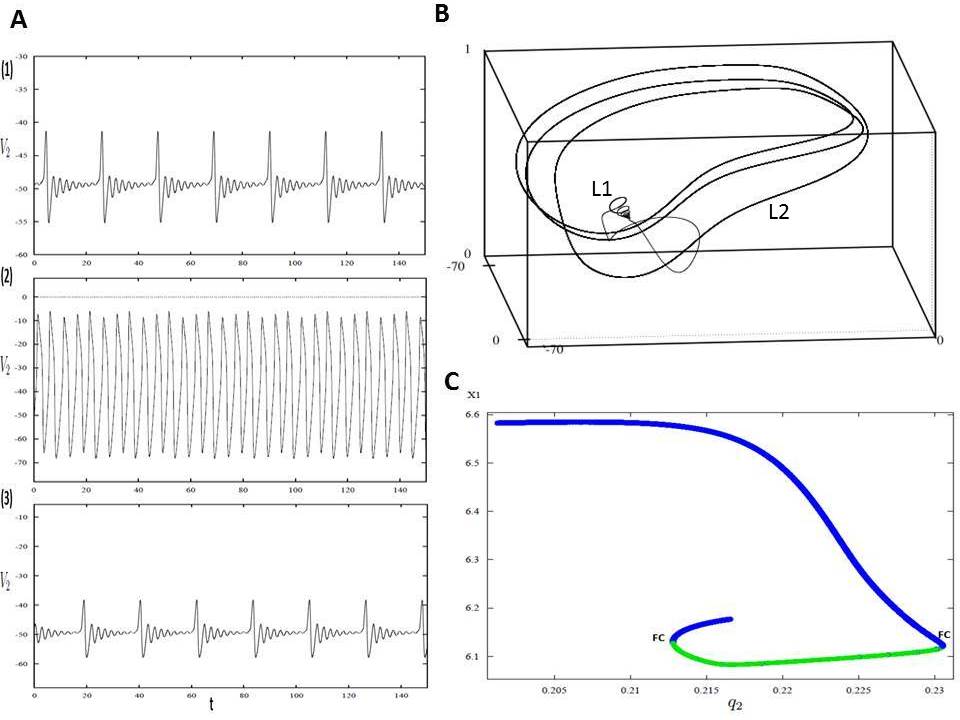} 
\caption{ A. The voltage time series of $"II"$ exhibits spiking mode for (2) $q_2=0.214$,   while at (1) $q_2=0.212$ and  (3) $q_2=0.24$,  $"II"$ exhibits subthreshold oscillations. B. For  $q_2=0.214$  the coupled system has two stable limit cycles,  $L_1$  and  $L_2$. The limit cycle $L_1$  corresponds to the subthreshold oscillations of $"II"$, while  $L_2$ corresponds to spiking mode in $"II"$.   C. Bifurcation diagram of $L_2$ which corresponds to the spiking mode in $"II"$, when $q_2$ is the bifurcation parameter and $X_1$ is the maximum of $V_1$ on the limit cycle. Here, the color green indicates stable limit cycle, the color blue  indicates saddle limit cycle and "FC" indicates the fold limit cycle bifurcation.}
\label{Figure10}
\end{center} 
\end{figure} 
By increasing the  coupling strength, subthreshold oscillations is observed in $"II"$.  More precisely, for $0.1019<q_2<0.23165$, the input signal to $"II"$ is not strong enough to make    $"II"$ fire an action potential. Hence,  $"II"$ exhibits subthreshold oscillations (Figure \ref{Figure9} B).

 A closer look shows that,  at $q_2=0.214$,  $"II"$  is in spiking mode, while at $q_2=0.212$ and $q_2=0.24$,  $"II"$ exhibits subthreshold oscillations (Figure \ref{Figure10} A).     
More precisely, consider the interval 
$(0.213,0.23165)$ in the parameter's space.   The following proposition holds for this interval:

\begin{proposition}
For the values of the coupling strength within the interval 
\\
$(0.213,0.23165)$  the resonator neuron is in tonic spiking mode, while for the values of the coupling strength outside which the resonator neuron exhibits subthreshold oscillations. 
\end{proposition}
\begin{proof}
For  $q_2=0.214$  the coupled system has two stable limit cycles  $L_1$  and  $L_2$ (Figure \ref{Figure10} B). The limit cycle $L_1$  corresponds to the subthreshold oscillations of $"II"$, while  $L_2$ corresponds to spiking mode in $"II"$.
The bifurcation diagram of $L_2$ has been depicted in Figure \ref{Figure10} C, where the bifurcation parameter is $q_2$.  As $q_2$ decreases, the stable limit cycle $L_2$ and a saddle one approach each other, coalesce at $q_2=0.213$ and  then disappear.  That is, at $q_2=0.213$ the stable limit cycle  $L_2$ undergoes fold limit cycle bifurcation, hence $L_2$ disappears. Then the state of the system tends to the other stable limit cycle $L_1$, hence subthreshold oscillations is observed in $"II"$. Now, again let $q_2=0.214$.  As $q_2$ increases, the stable limit cycle  $L_2$ tends to a saddle one, coalesce at $q_2=0.23165$  and then disappears. That is, by increasing $q_2$  the stable limit cycle $L_2$ undergoes fold limit cycle bifurcation, then it disappears. Hence, the  state of the system tends to the other stable limit cycle $L_1$, which corresponds to the subthreshold oscillations of $"II"$.
\end{proof}

The above proposition shows that, for $q_2=0.214$ the resonator neuron is in spiking mode, while for the coupling strength bigger or smaller than $q_2=0.214$, for instance $q_2=0.212$ and $q_2=0.24$, the resonator neuron exhibits subthreshold oscillations. That is, for $q_2=0.212$ and $q_2=0.24$,  the input signal to $"II"$ is not strong enough to make    $"II"$ fire an action potential, hence  $"II"$ exhibits subthreshold oscillations. 
\begin{remark}
This observation is because of that,   $"II"$, as a resonator neuron, may fire through fine tuning.
\end{remark}

\begin{remark}
By further increasing in the coupling strength, one can find some intervals in the parameter space, for the parameter values within which $"II"$ exhibits  bursting oscillations. More precisely, for the parameter values within these intervals the coupled system (\eqref{2}) has two stable limit cycles. One of them  corresponds  to  subthreshold oscillations of  $"II"$ and the other corresponds to  bursting oscillations of $"II"$.   
\end{remark} 

 \subsection{\textbf{Intermittent Oscillations}}
 Intermittency describes  a mechanism in many dynamical systems through which simple dynamics of the system becomes chaotic \cite{{diaz2001intermittency}}.
A conceptual formulation of intermittency  was proposed by Floris Takens in 1988 \cite{takens1988intermittancy}:

"A one-parameter family of  diffeomorphisms $\{\phi_{\mu}\}_{\mu}$ on a manifold has an intermittency bifurcation for $\mu=\mu_0$  at a compact invariant set $K$ if:

1. For every  $\mu< \mu_0$ the diffeomorphism $\phi_{\mu}$ has a compact invariant set $K_{\mu}$ converging to $K$ in Hausdorff sense when $\mu$ tends to $\mu_0$  from blow,

2. for $\mu>\mu_0$ close to $\mu_0$ there are no $\phi_{\mu}$--attracting sets near $K$, but the orbit of   Lebesgue almost every point in a neighbourhood of $K$ returns close to $K$ infinitely often."

The following theorem, which has been stated by Newhouse, Palis, Takens in 1983   \cite{newhouse1983bifurcations},  helps us to find intermittency bifurcations in the coupled system (\eqref{2}).

\begin{figure}[!ht]
\begin{center} 
\includegraphics[width=1\linewidth,height=1\linewidth]{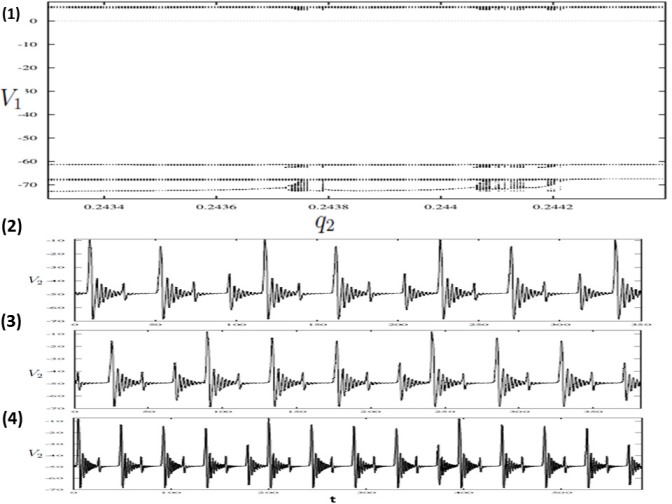} 
\caption{  By  increasing  $q_2$ the system undergoes  a sequence of  homoclinic bifurcations  to fold limit cycle,  through which   the coupled system (\eqref{2}) becomes chaotic. (1)   Poincar{\'e} map's bifurcation diagram  of  the system when  $q_2$ is the bifurcation parameter and the  Poincar{\'e} section is $V_2=-41.5$. (2),(3),(4). Voltage time series of $"II"$ for  (1) $q_2=0.243235$, (2) $q_2=0.243943$ and (3) $q_2=0.244165$.}
\label{Figure11}
\end{center} 
\end{figure}

\begin{theorem}
Suppose that a family of surface diffeomorphisms $\{\phi_{\mu}\}_{\mu}$ unfolds generically a non-transverse homoclinic bifurcation to fold limit cycle at $\mu=0$. Then there is a sequences $\nu_n\rightarrow0$  such that for every $\nu_n$, the diffeomorphism $\phi_{{\nu}_n}$ has a homoclinic tangency which is unfolded  generically by the family $\{\phi_{\mu}\}$.
\end{theorem}

For strong enough coupling strength, $"II"$ exhibits a new  oscillation pattern, which is  a combination of spiking oscillations and subthreshold oscillations.  In the following, it is examined that in the coupled system (\eqref{2}), transition  mechanism
between different observed oscillations of this pattern  gives rise to  intermittency.

Figure \ref{Figure11}(1) shows the Poincar{\'e} map{'}s bifurcation diagram of  the  coupled system (\eqref{2}). The bifurcation  diagram admits  that,   by  increasing  $q_2$ the system undergoes  a sequence of  homoclinic bifurcations  to fold limit cycle,  through which   the coupled system (\eqref{2}) becomes chaotic.  Each chaotic regime eventually ends with another homoclinic bifurcation  to fold limit cycle, then a new stable limit cycle appears. The new stable limit cycle  corresponds  to a  new combination of spiking oscillations and subthreshold oscillations. That is, the system has   an intermittency bifurcation for the  bifurcation values.

 Figure \ref{Figure11} (2),(3),(4)  show the voltage time series of $"II"$ which correspond  to  some of the observed intermittent oscillations. 

\begin{remark}
It is noticeable that by increasing the coupling strength, the length and the period of the stable limit cycle of the coupled system (\eqref{2}) increase.
\end{remark}

\subsection{\textbf{Synchronous Oscillations}}
For strong enough coupling strength,   the coupled oscillators show synchronous oscillations \cite{chow2000dynamics}. In neuroscience,  the cells synchronize when they oscillate with the same amplitude and
frequency without any phase differences.  

\begin{figure}[!ht]
\begin{center}
\includegraphics[width=.7\linewidth]{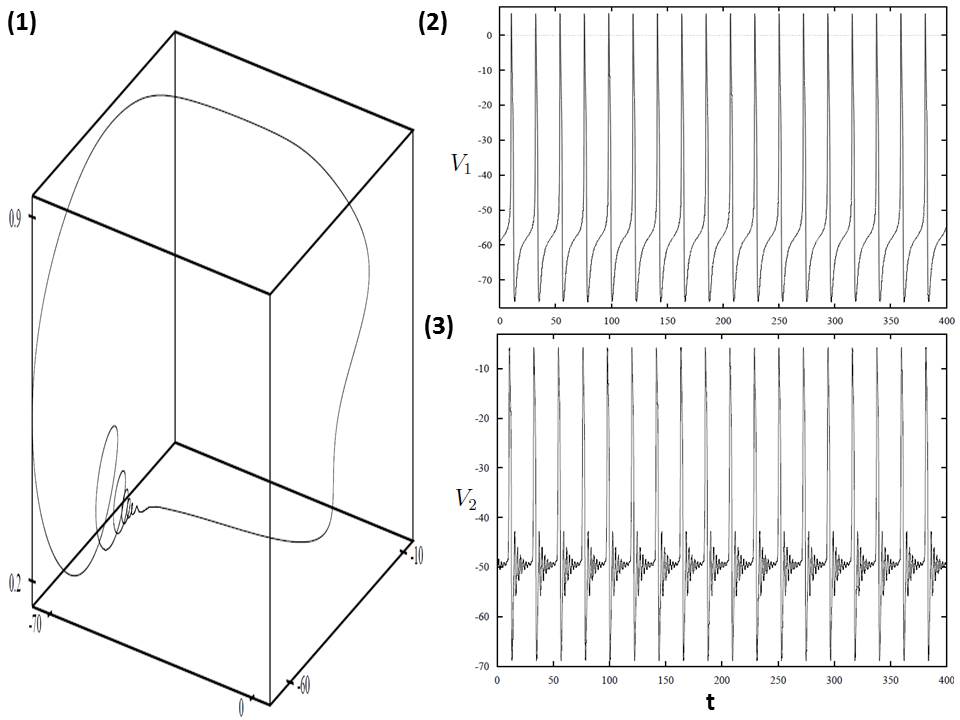} 
\caption{ For strong enough coupling strength,  synchronous oscillations  is observed in the coupled neurons $"I"$ and $"II"$.  (1) Three-dimensional image of the stable limit cycle which corresponds to the synchrony,  and  the corresponding  voltage time series of $"I"$ and $"II"$, for $q_2=0.35$.  }
\label{Figure12}
\end{center} 
\end{figure}  
For strong enough coupling strength, synchronous oscillations  is observed in the coupled neurons $"I"$ and $"II"$ (Figure \ref{Figure12}). A closer  inspection of the synchronous oscillations, reveals the following result:

\begin{figure}[!ht]
\begin{center}
\includegraphics[width=.8\linewidth]{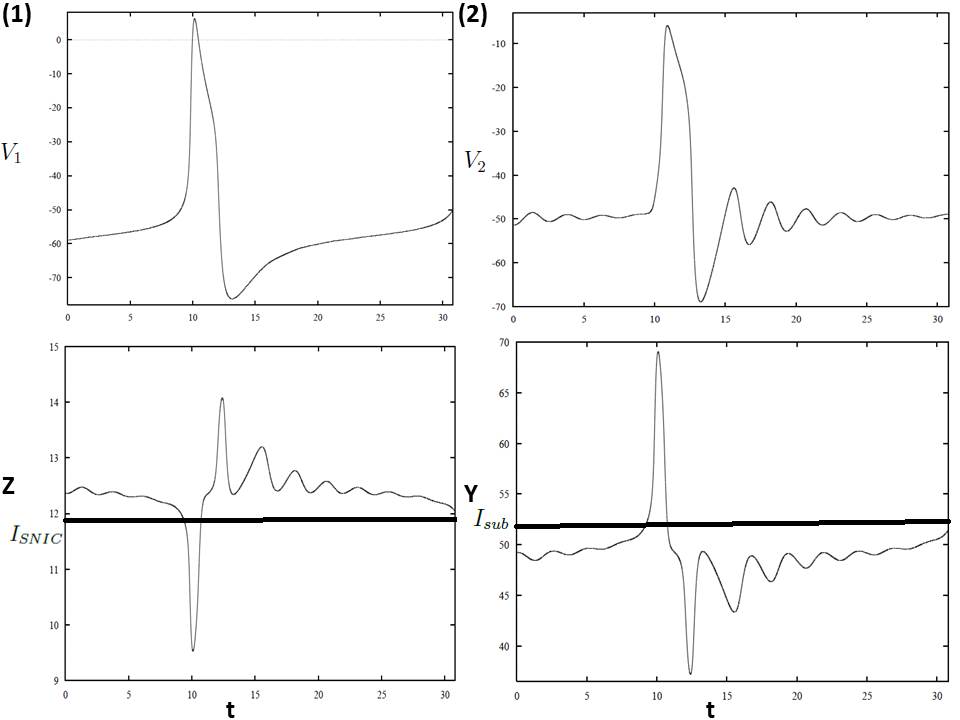} 
\caption{ The input signals through the input current and   the coupling  strength to $"I"$ and $"II"$,   which are obtained by $Y=I_1+q_1(V_2-V_1)$ and $Z=I_2+q_2(V_1-V_2)$ respectively, For $q_2=0.35$. There exist some time intervals within which the input signal to $"II"$ is such that   $"II"$ fires single spike and then comes back to its silent phase.   Hence,   $"II"$ exhibits 1-Bursting oscillations.
 Moreover, the input signals to $"I"$ is not enough to make $"I"$ become silent.  Hence, periodic spiking of $"I"$ persists.  }
\label{Figure13}
\end{center} 
\end{figure}  

\begin{proposition}
In the coupled system  (\eqref{2}), the final synchronization of $"I"$ and $"II"$ actually corresponds to the synchronization of tonic spiking oscillations of $"I"$ and 1-bursting oscillations of $"II"$.
\end{proposition}

\begin{proof}
Let $I_{SNIC}:=11.99$ and $I_{sub}:=51.9$. As mentioned in the section $2.1$, for $m'_{1/2}=-30$,  
 at $I=I_{SNIC}$ the neuron (i.e. $"I"$)  undergoes saddle-node bifurcation on invariant circle, hence $"I"$ initiates spiking oscillations. Moreover, for $m'_{1/2}=-45$,
at $I=I_{sub}$  the neuron (i.e. $"II"$) undergoes subcritical Hopf bifurcation, then $"II"$ initiates spiking oscillations.
As mentioned previously, for $q_2=0.35$, $"I"$ and $"II"$ exhibit synchronous oscillations (Figure \ref{Figure12}). 
  Suppose that  $Y$ and $Z$ are input signals through the input current and   the coupling  strength to $"I"$ and $"II"$,   which are obtained by $Y=I_1+q_1(V_2-V_1)$ and $Z=I_2+q_2(V_1-V_2)$ respectively. As demonstrated by Figure \ref{Figure13} (2), there exist some time intervals within which the input signal to $"II"$ is less than the corresponding bifurcation value $I_{sub}$, hence damping oscillations  of $"II"$ are observed. By passing through the time interval, the input signal to $"II"$ is such that, $"II"$ fires single spike and then comes back to its silent phase.   Hence,   $"II"$ exhibits 1-Bursting oscillations.
 Moreover, the input signals to $"I"$ is not enough to make $"I"$ become silent (Figure \ref{Figure13} (1)).  Hence, periodic spiking of $"I"$ persists. 
  Therefore, when $q_2=0.35$,
 the synchrony is actually the synchronization of periodic spiking of $"I"$  and 1-bursting oscillations of $"II"$.
\end{proof}

So far,     all the observed oscillation patterns of the coupled neurons $"I"$ and $"II"$ have been introduced. In the next section,  transition mechanisms between different patterns  are examined.

\section{Transition Mechanisms Between Different oscillation patterns}

 \begin{figure}[h!]
\begin{center}
\includegraphics[width=.7\linewidth,height=.6\linewidth]{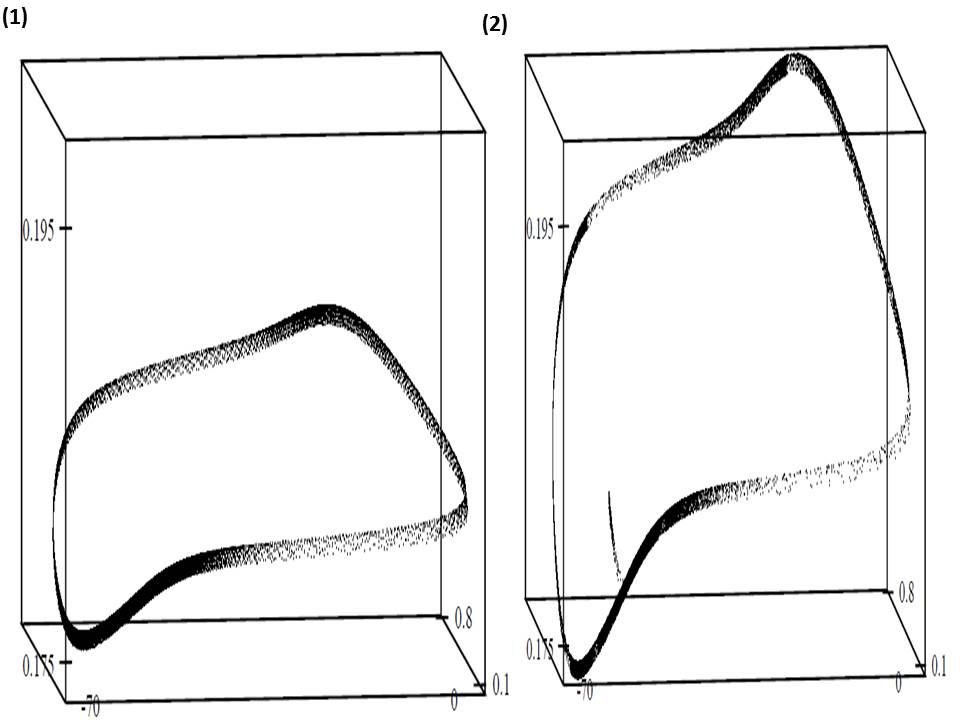} 
\caption[]{  By increasing the coupling strength, the smooth torus tends to lose its smoothness.  The image of the  Poincar\'e map on  $(1)$ the smooth torus for $q_2=0.04$ and $(2)$  the non-smooth torus for  $q_2=0.8616$, where 
$V_2=-50$
is the corresponding Poincar\'e section on the torus.  }
\label{Figure14}
\end{center} 
\end{figure} 

 In the previous  section it was stated that, in the coupled system (\eqref{2}) increasing the coupling strength results in   different oscillation patterns of $"II"$. In this  section,   transition mechanisms between different oscillation patterns  are investigated.
\subsection{\textbf{Transition From Phase-Locking Oscillations to TS/MMOs}}
 As mentioned in the section $3.1$,  when the solution of the coupled system (\eqref{2}) lies on the smooth two-dimensional torus, then the solution corresponds to the phase-locking oscillations of $"I"$ and $"II"$. Since the torus is smooth,  all the Poincar\'e maps that correspond to suitable sections, are smooth (Figure \ref{Figure14} (1)). By increasing the coupling strength, the smooth torus tends to lose its smoothness. 
 Figure \ref{Figure14} (2) shows the image of the non-smooth Poincar\'e map for $q_2=0.08616$, hence the torus is non-smooth. Now, by using  "Affraimovich-Shilnikov break-down theorem" \cite{afraimovich1982bifurcation} one can conclude that:
 
 \begin{figure}[h!]
\begin{center}
\includegraphics[width=1\linewidth]{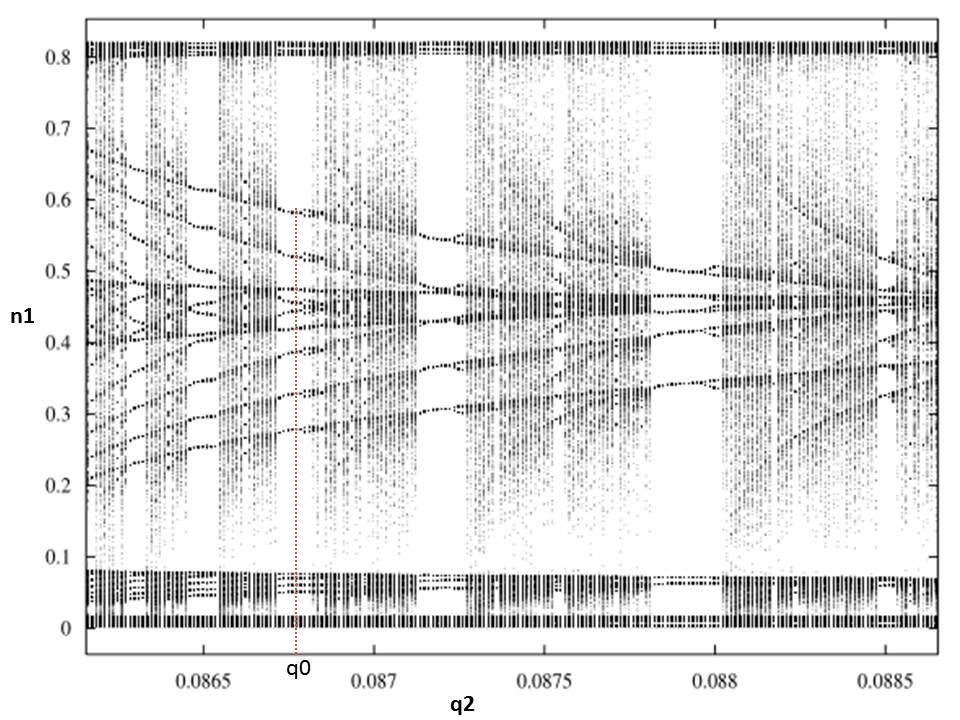} 
\caption[]{   Bifurcation diagram of the Poincar\'e map of the coupled  system (\eqref{2}) corresponding to the section $V_2=-50$, where $q_2$ is the bifurcation parameter.    At $q_0:=q_2=0.086785$ the stable limit cycle undergoes period-doubling bifurcation. The occurrence of this bifurcation shows that the invariant torus no longer exists.  }
\label{Figure15}
\end{center} 
\end{figure} 

\begin{proposition}
By increasing $q_2$ the stable torus breaks down. More precisely, first the torus loses its smoothness and then, it breaks down through homoclinic bifurcation to a saddle-node limit cycle.
\end{proposition}

\begin{proof}
By increasing $q_2$ the stable torus loses its smoothness  (Figure \ref{Figure14} (2)). Figure \ref{Figure15}  shows the bifurcation diagram of the Poincar\'e map of the coupled system (\eqref{2}). By further increasing in the coupling strength $q_2$,  at  $q_0:=q_2=0.086785$ the coupled system has a   cycle 
such that one of whose  flouquet multipliers is equal to $-1$, then the limit cycle becomes unstable. That is, at $q_0$ the stable limit cycle undergoes period-doubling bifurcation (Figure \ref{Figure15}). The occurrence of this bifurcation shows that the invariant torus no longer exists. 

By using "Affraimovich-shilnikov break-down theorem" it is concluded that, there exists some $q^* \in (0.086153,0.086785)$ for which the torus breaks down through homoclinic bifurcation to saddle-node  cycle. More precisely, at $q_2=q^*$ the torus becomes global unstable set of the saddle-node cycle and  then the torus breaks down.
\end{proof}

After the torus destruction, a stable limit cycle $M$ appears.  It is noticeable that, $M$ corresponds   to tonic spiking oscillations of $"I"$ and  mixed mode oscillations   of $"II"$ (Figure \ref{Figure4}), hence the torus destruction results in a new   oscillation pattern, i.e. TS/MMOs.  

 \begin{figure}[h!]
\begin{center}
\includegraphics[width=1\linewidth,height=.6\linewidth]{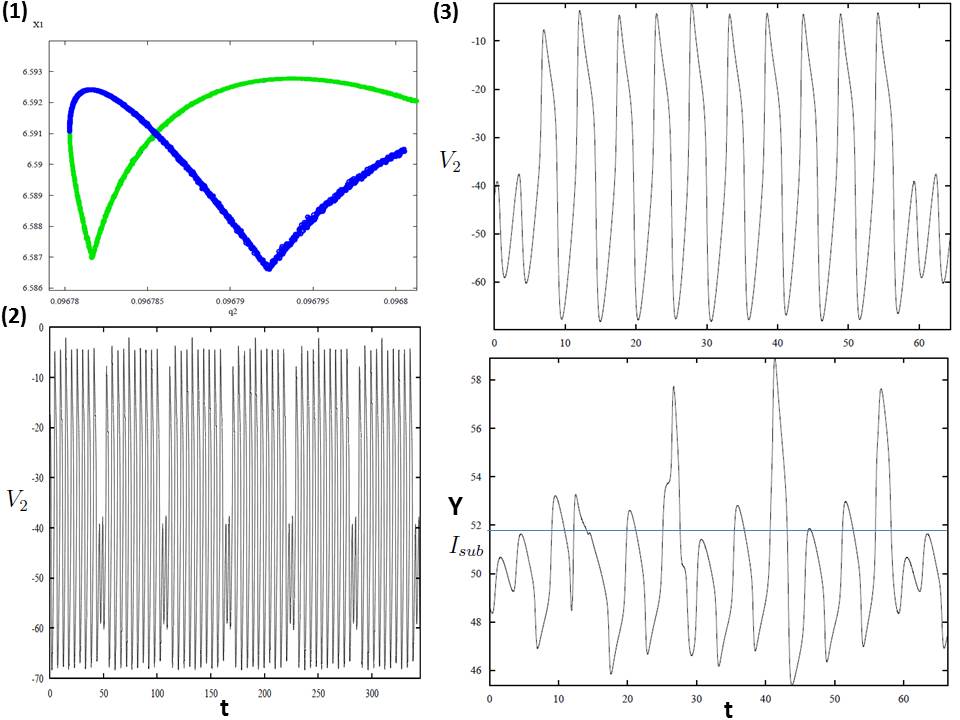} 
\caption[]{  By increasing the coupling strength $q_2$, at $q_2 = 0.096779$ the coupled system (\eqref{2}) undergoes  fold  limit cycle bifurcation, then a new stable limit cycle appears. This limit cycle corresponds to the bursting oscillations of $"II"$. (1) The bifurcation diagram of the coupled system (\eqref{2}), when  the coupling strength $q_2$ is the bifurcation parameter and   $X_1$ is the maximum of $V_1$ on the limit cycle. (2) The voltage time series of "II", $V_2$,  corresponding to the new stable limit cycle for $q_2 = 0.09678$. (3) The corresponding input signal, $Y$,  to $"II"$  
 through the input current and   the coupling  strength, i.e. $Y=I_2+q_2(V_1-V_2)$, during the active phase of $"II"$. There exist some time intervals within which the input signal to $"II"$ is below the corresponding bifurcation value $I_{sub}$, hence in each interval $"II"$ is silent.}
\label{Figure16}
\end{center} 
\end{figure}

 \subsection{\textbf{Transition From TS/MMOs  to TS/Bursting}}

As mentioned previously in the section 3.2, for strong enough coupling strength the resonator neuron exhibits mixed mode oscillations. 
\begin{proposition} 
By  increasing  the coupling strength, one can find an interval in the parameter space for the parameter values within which $"II"$ may exhibit bursting oscillations. 
 \end{proposition}
\begin{proof}
 The bifurcation diagram of the system has been depicted in Figure \ref{Figure16} (1), when $q_2$ is the bifurcation parameter. As shown by this diagram, at $q_2 = 0.096779$ the coupled system (\eqref{2}) has a saddle-node  cycle, then by increasing the coupling strength a stable limit cycle appears. That is, at $q_2 = 0.096779$ the system undergoes fold limit cycle bifurcation and then a stable limit cycle appears.
 Figure \ref{Figure16} (2) 
shows the voltage time series of $"II"$, which corresponds to this stable limit cycle for $q_2 = 0.09678$. In the following, it is investigated that this stable limit cycle corresponds to the bursting oscillations of $"II"$.

 Let $q_2 = 0.09678$ and $I_{sub} := 51.9$. As mentioned in the section 2.1, at $I = I_{sub}$, $"II"$ undergoes subcritical Hopf bifurcation, then $"II"$ initiates spiking oscillations. Let $Y$ be the input signal to $"II"$ through the input current and the coupling strength, which is obtained by $Y = I_2 +q_2(V_1-V_2)$. As shown by Figure \ref{Figure16} (3),

 \begin{figure}[h!]
\begin{center}
\includegraphics[width=1\linewidth,height=.65\linewidth]{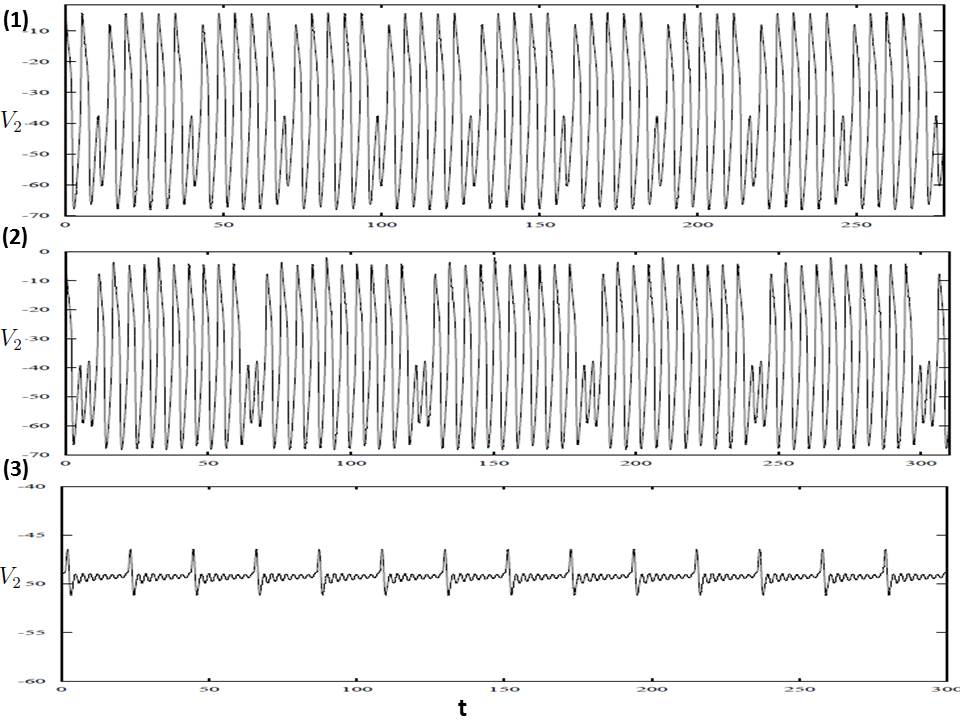} 
\caption[]{For $q_2=0.096783$, $"II"$ exhibits (1)  mixed mode oscillations when $x_0=(-62.2114,0.0027,-5.4806,0.6079)$, (2)   bursting oscillations when
$x_0=(-55.1947,0.0051,-4.5256,0.6247)$ and (3) subthreshold oscillations when $x_0=(-62.2114,0.0027,-5.4806,0.6079)$. Here  $x_0$ is the initial condition. }
\label{Figure17}
\end{center} 
\end{figure} 
there exist some time intervals within which, the input signal to $"II"$ is below the corresponding bifurcation value $I_{sub}$, hence in each interval $"II"$ is silent. Going through each of these time intervals, the input signal to $"II"$ becomes bigger than the corresponding bifurcation value $I_{sub}$, therefore $"II"$ starts periodic spiking. The periodic spiking continues until that at the next time interval $"II"$ becomes silent again. In the other words, $"II"$ shows bursting oscillations. Since the input signal to $"II"$ is periodic, periodic bursting oscillations is observed in $"II"$.
\end{proof}

 \begin{remark} 
It is noticeable that, for $q_2 > 0.096779$ the coupled system (\eqref{2}) is multistable. That is, the system (\eqref{2}) has three stable limit cycles, which correspond to the mixed mode oscillations of $"II"$  (Figure \ref{Figure17} (1)), the bursting oscillations of $"II"$ (Figure \ref{Figure17} (2)) and the subthreshold oscillations of $"II"$ (Figure \ref{Figure17} (3)).
\end{remark}

 \begin{figure}[h!]
\begin{center}
\includegraphics[width=1\linewidth,height=.65\linewidth]{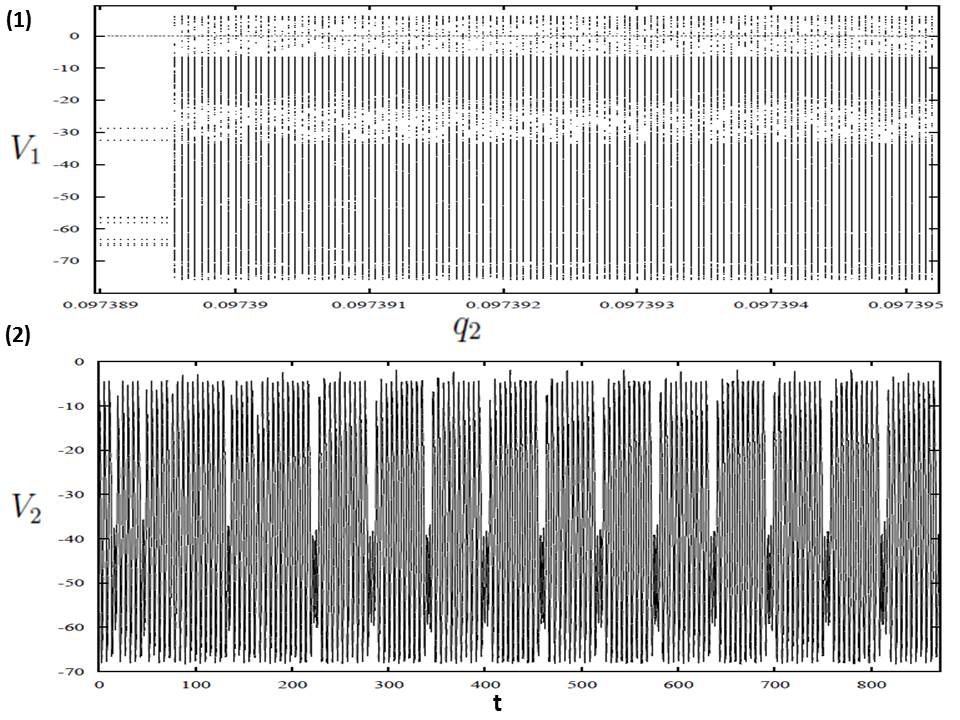} 
\caption[]{   At $q_2 = 0.0973895$, the stable limit cycle which corresponds to the MMOs, undergoes fold limit cycle bifurcation, then it disappears. Then for suitable initial conditions, $"II"$ exhibits bursting oscillations. $(1)$ Bifurcation diagram of the Poincar\'e map of the system corresponding to the section $V_2 =-30$, when $q_2$ is the bifurcation parameter. $(2)$ The voltage time series of $"II"$ corresponding to the bursting oscillation   
for $q_2=0.0.09739$.  }
\label{Figure18}
\end{center}  
\end{figure} 
 
In the following, it is investigated that how the mixed mode oscillations end up.
\begin{proposition}
By further increasing in the coupling strength, the mixed mode oscillations of $"II"$ end up through fold limit cycle bifurcation.
\end{proposition}
\begin{proof}
Figure \ref{Figure18}(1) shows the bifurcation diagram of the Poincar\'e map of the system corresponding to the section $V_2 = -30$, when $q_2$ is the bifurcation parameter. As depicted by this diagram, at $q_2 = 0.0973895$ the stable limit cycle which corresponds to the MMOs, undergoes fold limit cycle bifurcation, then it disappears. Hence, depending on the initial condition, $"II"$ may exhibit bursting oscillations  (Figure \ref{Figure18} (2)) or subthreshold oscillations.
\end{proof}

 \begin{figure}[h!]
\begin{center}
\includegraphics[width=1\linewidth,height=.6\linewidth]{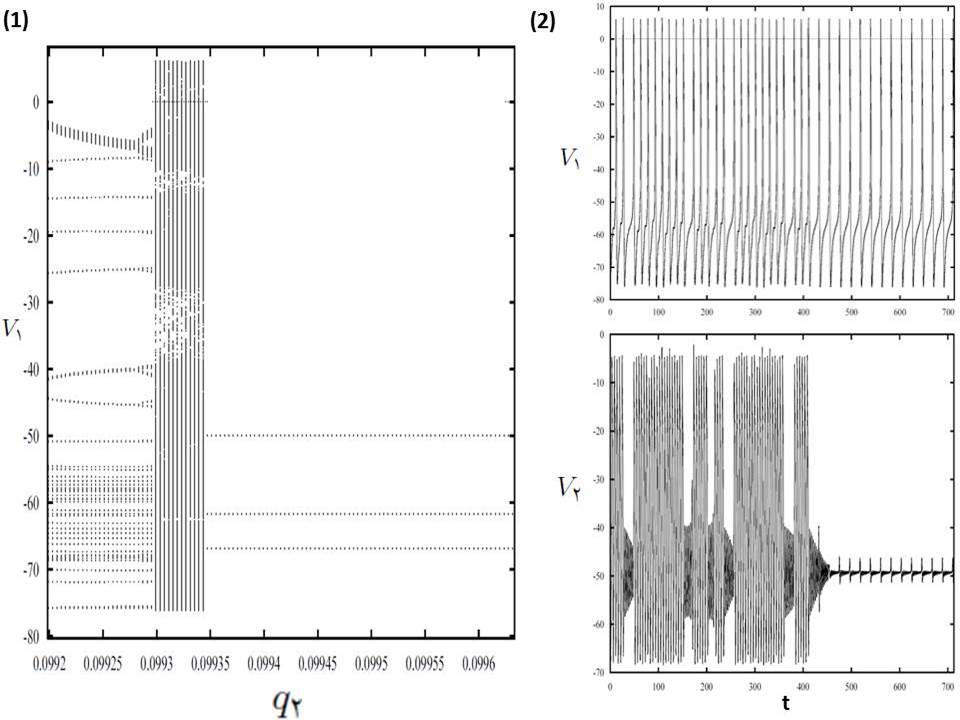} 
\caption[]{ For 
$q_2>0.09935$
the stable limit cycle which corresponds to the subthreshold  oscillations of $"II"$ is globally stable. That is, starting from arbitrary initial conditions, the solution eventually tends to the limit cycle, which corresponds to the subthreshold oscillations of   $"II"$. $(1)$ Bifurcation diagram of the Poincar\'e map of the system corresponding to the section $V_2=-40$, when $q_2$ is the bifurcation parameter, and 
$(2)$ the voltage time series of $"I"$
 and
  $"II"$  
for $q_2=0.1019$.  }
\label{Figure19}
\end{center}  
\end{figure} 

 \subsection{\textbf{Transition From TS/Bursting to TS/Sub Oscillations}}

In the section 3.3, it has been stated that, for strong enough coupling strength the coupled system is bistable. That is, the system has two stable limit cycles $U_1$ and $U_2$. The limit cycle $U_1$ corresponds to the bursting oscillations of $"II"$, and $U_2$ corresponds to the subthreshold oscillations of $"II"$ (Figure \ref{Figure7} B). 

 \begin{figure}[!ht]
\begin{center}
\includegraphics[width=1\linewidth,height=.6\linewidth]{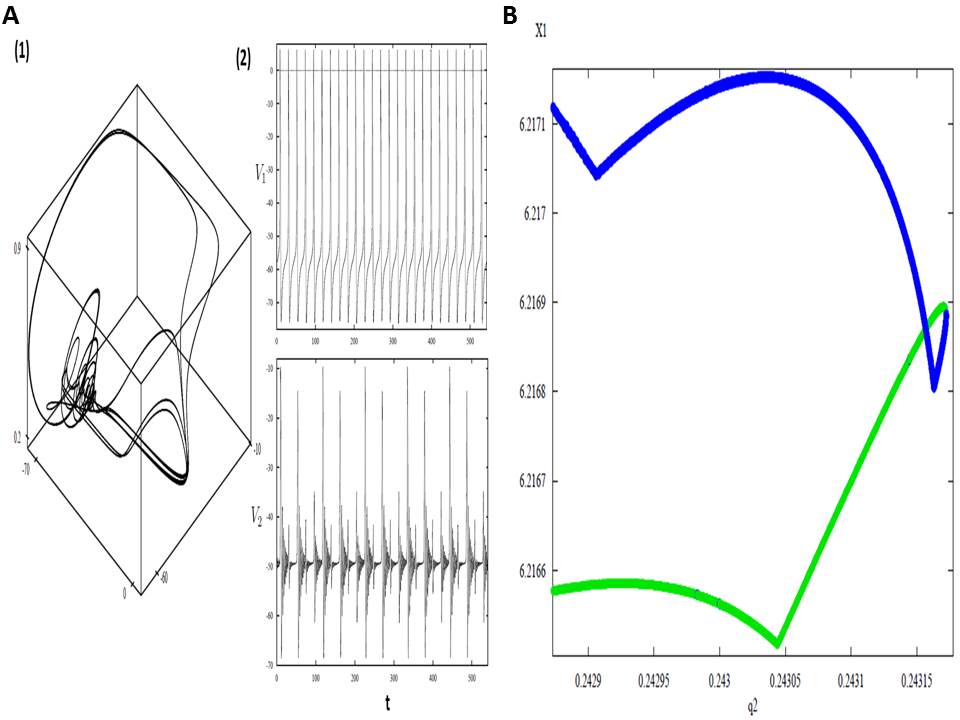} 
\caption[]{A. (1),(2). For $q_2=0.24317$, the coupled system (\eqref{2}) has a stable limit cycle which  corresponds to a combination of  bursting oscillations and subthreshold oscillations of $"II"$. (1) Three-dimensional image of the stable limit cycle  and  $(2)$ the corresponding voltage time series of $"I"$ and $"II"$.  B. Bifurcation diagram of the stable limit cycle which corresponds  to subthreshold oscillations of $"II"$, when $q_2$ is the bifurcation parameter and $X_1$ is the maximum value of $V_1$ on the  limit cycle. By increasing the coupling strength the stable limit cycle undergoes the fold limit cycle bifurcation, then it disappears.   }
\label{Figure20}
\end{center} 
\end{figure}

\begin{proposition}
 By increasing the coupling strength bursting oscillations of $"II"$ eventually end up. 
\end{proposition}
\begin{proof}
 Figure \ref{Figure19} (1) shows the bifurcation diagram of the Poincar\'e map of $U_1$. As depicted by the bifurcation diagram, by further increasing in the coupling strength the bursting oscillations end up. Then the state of the system tends to the other stable limit cycle $U_2$, which corresponds to the subthreshold oscillations of $"II"$. Hence, for $q_2 > 0.09935$ the stable limit cycle, which corresponds to tonic spiking of $"I"$ and subthreshold oscillations of $"II"$, is globally stable.  (Figure  \ref{Figure19} (2))
\end{proof}

  \subsection{\textbf{Transition From TS/Sub Oscillations to Intermittent Oscillations}}
 For strong enough coupling strength, the system has two stable limit cycles. One of them corresponds to the subthreshold oscillations, and the other corresponds to a combination of bursting oscillations and subthreshold oscillations. By increasing the coupling strength, the bursting oscillations eventually end up, then the resonator neuron exhibit subthreshold oscillations. In the following, this transition will be investigated precisely.

 \begin{proposition}
  By increasing the coupling strength,  the stable limit cycle, which corresponds to the subthreshold oscillations of $"II"$,  annihilates through  fold limit cycle bifurcation.
  \end{proposition}
\begin{figure}[ht]
\begin{center}
\includegraphics[width=1\linewidth,height=.5\linewidth]{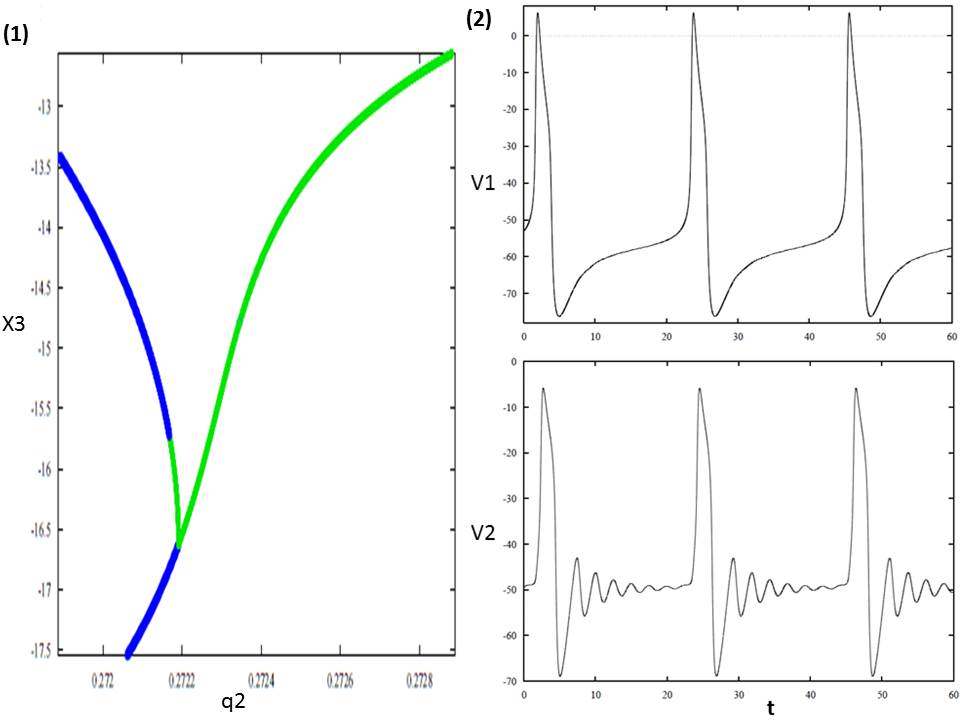} 
\caption[]{$(1)$  Bifurcation diagram of the stable  limit cycle,  which corresponds to  intermittent oscillations of $"II"$, when $q_{2}$ is the bifurcation parameter and $X_3$ is the maximum value of   $V_2$ on the cycle.    By increasing $q_2$, at  $q_2=0.2722$ the stable limit cycle undergoes homoclinic bifurcation to a saddle-node cycle. Then, it disappears and  a new stable limit cycle, which corresponds to the synchrony, appears.
  $(2)$ The  voltage time series of the coupled neurons corresponding to the synchronous oscillations. }
\label{Figure21}
\end{center} 
\end{figure}   
   \begin{proof}
    For $q_2=0.24317$, the coupled system (\eqref{2}) has a stable limit cycle, which corresponds to a combination of bursting oscillations and subthreshold oscillations  of $"II"$ (Figure \ref{Figure20} A(1),(2)).  The  bifurcation diagram of the stable limit cycle, which corresponds to the subthreshold oscillations of $"II"$, has been depicted in  Figure \ref{Figure20} B,
when $q_2$ is the bifurcation parameter. As demonstrated by this diagram, by increasing the coupling strength this stable limit cycle and a saddle one approach each other, they collide at $q_2 = 0.24317$ and then annihilate each other. That is, at  $q_2 = 0.24317$ the stable limit cycle undergoes the fold limit cycle bifurcation, then it disappears. After this bifurcation, the other stable limit cycle, which corresponds to a combination of subthreshold oscillations and bursting oscillations of $"II"$, is globally stable.  
 \end{proof}

 \subsection{\textbf{Transition From Intermittent Oscillations to Synchronous Oscillations}}

As previously explained in the section 3.6, for strong enough coupling strengths, the coupled neurons $"I"$ and $"II"$ exhibit synchronous oscillations. The following proposition, talks about the transition mechanism between intermittent oscillations and synchronous oscillations of the neurons.

\begin{proposition}
  The intermittent oscillations of $"II"$ eventually end up through homoclinic bifurcation to a saddle-node cycle. Then, the coupled neurons exhibit synchronous oscillations.
\end{proposition}   

\begin{proof} 
Figure \ref{Figure21} (1) shows the bifurcation diagram of  the stable limit cycle, which corresponds to  intermittent oscillations of $"II"$, where   $q_2$ is the bifurcation parameter. As demonstrated by this diagram,  by increasing the coupling strength the stable limit cycle and a saddle one  approach each other and collide at $q_2=0.2722$.
By further increasing in the coupling strength,  the saddle-node cycle annihilates, then a new stable limit cycle appears. That is, by increasing $q_2$, at  $q_2=0.2722$ the stable limit cycle undergoes homoclinic bifurcation to a saddle-node cycle. Then, it disappears and a new stable limit cycle appears. This limit cycle corresponds to the synchronization of $"I"$ and $"II"$ (Figure \ref{Figure21} (2)).
\end{proof}

\section*{\textbf{Discussion}}
A great deal of research has been devoted on the investigation of neural dynamics in coupled neurons. One of  the  most important questions is that, in a network of  neurons how the differences  in the   dynamics of the neurons  affects the dynamics of the network. In \cite{razvan2020} the dynamics of two coupled integrator neurons of different types of excitability through the gap junction has been investigated. The results of \cite{razvan2020} shows that increasing the coupling strength between the coupled neurons determines a rich dynamic behaviour, such as emergent bursting oscillations of the coupled neurons and burst synchronization of the coupled neurons. 

One important question is that how differences
in intrinsic characteristics of the coupled neuron’s dynamics, e.g. resonance
or integratory dynamics, affects the network's dynamics, while the neurons are of different types of excitability.  In this manuscript, the answer to this question  has been  investigated,  when a silent resonator neuron is coupled with a spiking integrator neuron through a gap junction. Moreover, none of the neurons exhibit mixed mode oscillations and bursting oscillations intrinsically.  Then, by using the dynamical systems theory (e.g. the bifurcation theory), it has been  examined that how increasing the coupling strength affects the dynamics of the neurons, when one of the coupling strength is fixed and the other varies.

 By increasing the coupling strength, different oscillation patterns, such as mixed mode oscillations and bursting oscillations, have been   observed in the resonator neuron, while the integrator neuron is in tonic spiking mode. 
By increasing the coupling strength,  multi-stability, as one of the most important features of the coupled system, has been observed  in the system. That is, for strong enough coupling strength the coupled system  has three stable limit cycles,  which correspond to the mixed mode oscillations,  bursting oscillations   and   subthreshold oscillations of the resonator neuron. In the coupled system, the bursting oscillations  have  a feature  that distinguishes them from  most of the other  observed bursting oscillations. In   the bursting oscillations, usually   the stable limit  cycle, which corresponds to the subthreshold oscillations, loses its stability, though in the coupled system the stability of the limit cycle persists. 
  the final synchronization of the neurons actually corresponds to the synchronization of tonic spiking oscillations of the integrator neuron and 1-bursting oscillations of the resonator neuron.
 The achievements  of this manuscript also confirm that,   the resonator neurons may fire through the fine tuning. More precisely, some  interval in the parameter space has been introduced  such that for the  values of the coupling strength  within this interval the resonator neuron is in spiking mode, while for the values of the coupling strength outside of which   the resonator neuron exhibits the subthreshold oscillations.  
 
In conclusion, given the rather complex oscillation patterns of
two coupled neurons, when the uncoupled neurons have no intrinsic mixed mode oscillations and bursting mode behaviour, one expects a much more complex behaviour when a
network of more than two neurons are coupled together, while at least one of the  coupled neurons is an integrator neuron and one of which is a resonator neuron and the coupled neurons are of different types of excitability.

\section*{Acknowledgement}
The authors  would  like to especially thank  Dr. Abdolhosein Abbasian  for his  helpful comments.

\bibliographystyle{alpha}

\end{document}